\documentclass[11pt,a4paper]{article}

\usepackage{amsmath,amssymb,amsthm,float}
\usepackage{graphicx}

\usepackage{amssymb,amsmath,amsthm,amscd,amsfonts}
\usepackage{amsmath,amssymb,amsthm,amscd,ascmac,mathrsfs,bm,type1cm,multirow,array,enumerate,mathrsfs,tabularx}
\usepackage{graphics}
\usepackage{tikz}

\numberwithin{equation}{section}
\numberwithin{figure}{section}

\newtheorem{theorem}{Theorem}[section]
\newtheorem{lemma}[theorem]{Lemma}
\newtheorem{proposition}[theorem]{Proposition}

\newcommand{\ve}{\bm{e}} 

\newcommand{\vX}{\bm{X}}

\begin{document}

\title{\textbf{Translators invariant under hyperpolar actions}}
\author{Tomoki Fujii and Naoyuki Koike}
\date{}
\maketitle

\footnote[0]{
\hspace{0.6cm}Key words and phrases.  Translator, Hyperpolar action, Mean curvature flow.
}

\begin{abstract}
In this paper, we consider translators (for the mean curvature flow) given by a graph of a function on a symmetric space $G/K$ of compact type which is invariant under 
a hyperpolar action on $G/K$.  
First, in the case of $G/K=SO(n+1)/SO(n)$, $SU(n+1)/S(U(1)\times U(n))$, $Sp(n+1)/(Sp(1)\times Sp(n))$ or $F_4/{\rm Spin}(9)$, we classify the shapes of translators in $G/K\times\mathbb R$ 
given by the graphs of functions on $G/K$ which are invariant under the isotropy action $K\curvearrowright G/K$.  
Next, in the case where $G/K$ is of higher rank, we investigate translators in $G/K\times\mathbb R$ given by the graphs of functions on $G/K$ which are invariant under 
a hyperpolar action $H\curvearrowright G/K$ of cohomogeneity two.  
\end{abstract}

\section{Introduction}
Let $(N,g)$ be a $n$-dimensional Riemannian manifold and $u:M\to\mathbb{R}$ be a function on a domain $M\subset N$.  
The {\it graph embedding} $f$ for $u$ is defined as the embedding of $M$ into the product Riemannian manifold $N\times\mathbb{R}$ given by $f(x):=(x,u(x))\,\,\,(x\in M)$.  
For the simplicity, denote by $\Gamma$ the graph $f(M)$ of $u$.  
If a $C^{\infty}$-family $\{f_t\}_{t\in I}$ of $C^{\infty}$-immersions of $M$ into $N\times\mathbb{R}$ ($I$ is an open interval including $0$) satisfies 
\begin{equation}\label{eq:mcf}
\begin{cases}
\displaystyle\left(\frac{\partial f_t}{\partial t}\right)^{\bot_{f_t}}=H_t\\
f_0=f,
\end{cases}
\end{equation}
the family $\{M_t\}_{t\in I}$ of the images $M_t:=f_t(M)$ is called the {\it mean curvature flow starting from} $\Gamma$, where 
$H_t$ is the mean curvature vector field of $f_t$ and $(\bullet)^{\bot_{f_t}}$ is the normal component of $(\bullet)$ with respect to $f_t$.  
According to Hungerb\"{u}hler and Smoczyk\cite{HS}, we define a soliton for the mean curvature flow as follows.  
Let $\vX$ be a Killing vector field on $N\times\mathbb{R}$ and $\{\phi_t\}_{t\in\mathbb{R}}$ be the one-parameter transformation group associated to $\vX$.  
If $\{f_t\}_{t\in I}$ satisfies 
\begin{equation}
\left(\frac{\partial(\phi_t^{-1}\circ f_t)}{\partial t}\right)^{\bot_{\widetilde{f}_t}}=0,\label{eq:mcf-soliton}
\end{equation}
then $\Gamma$ is called a \textit{soliton for the mean curvature flow with respect to $\vX$}.  
In the sequel, we call such a soliton a \textit{$\vX$-soliton} simply.  
In particular, when $X=(0,1)\in T(N\times\mathbb{R})(=TN\oplus T\mathbb{R})$, we call the $\vX$-soliton a \textit{translator}.  

The translator in $\mathbb{R}^{n+1}(=\mathbb R^n\times\mathbb R)$ has been studied by many researchers.  
Minimal submanifolds is special case of translators.
Therefore, translators can be regarded as the generalization of minimal submanifolds.
Hence, Bao and Shi \cite{BS} showed a Bernstein-type theorem for complete translators in case of codimension one.
Also, by Kunikawa \cite{Ku}, a Bernstein-type theorem for complete translators with flat normal bundle in case of higher codimension was shown.
Clutterbuck, Schn\"{u}rer and Schulze \cite{CSS} showed the existence of the complete rotationally symmetric graphical translator, which is called a {\it bowl soliton}.  
Here we note that Altschuler and Wu \cite{AW} had earlier showed its existence in the case of $n=2$
Also, they \cite{CSS} investigated a certain kind of stability for the bowl soliton.  
Denote by $r$ the distance function from the base point of $\mathbb{R}^n$.  Assume that the graph $\Gamma$ of a function $u$ on $M$ gives a bowl soliton.  
Then, since the bowl soliton is a rotationally symmetric entire graph, $u$ is described as $u=V\circ r$ for some function $V:[0,\infty)\to\mathbb R$.  
Furthermore, they showed that the function $V$ is a solution of some ODE and that $V(r)$ has the following asymptotic expansion as $r\to\infty$:
$$
\frac{r^2}{2(n-1)}-\log{r}+O(r^{-1}),
$$
Wang \cite{W} showd that when $n=2$, the bowl soliton is the only convex translator given as an entire graph.  
Further, Spruck and Xiao \cite{SX} showed that the bowl soliton  with $n=2$ is the only complete translator which is an entire graph.
The level sets of the function $r$ gives an isoparametric foliation on $\mathbb R^n$.  
By noticing this fact, Fujii \cite{F} classified the shape of the translator given by a function on a unit sphere which is constant along each leaf of isoparametric foliation.
Each leaf is an isoparametric hypersurface and the numaber $k$ of its distinct principal curvatures is equal to $1$, $2$, $3$, $4$ or $6$ (see \cite{M}).
In cases of $k=1$, $2$, $3$, Cartan\cite{C} classified the isoparametric hypersurfaces.
The hypersurfaces are $S^{n-1}\subset S^n$ in case of $k=1$, $S^k\times S^{n-k-1}\subset S^n$ in case of $k=2$ and the tubes over the Veronese surfaces $\mathbb{R}P^2\subset S^4$, $\mathbb{C}P^2\subset S^7$, $\mathbb{Q}P^2\subset S^{13}$, $\mathbb{O}P^2\subset S^{25}$ (i.e., the principal orbits of the isotropy representations of the rank two symmetric spaces $SU(3)/SO(3)$, $(SU(3)\times SU(3))/SU(3)$, $SU(6)/Sp(3)$, $E_6/F_4$) in case of $k=3$.
These hypersurfaces are homogeneous.
In case of $k=6$, the hypersurfaces are homogeneous by the result of Dorfmeister and Neher\cite{DN} and Miyaoka\cite{Mi}.
The hypersurfaces are the principal orbits of the isotropy representations of $(G_2\times G_2)/G_2$, $G_2/SO(4)$.
In case of $k=4$, Ozeki and Takeuchi\cite{OT1,OT2} found that non-homogeneous isoparametric hypersurfaces are constructed as the regular level sets of the restrictions of the Cartin-M\"{u}nzner polynomial functions to the sphere.
When the isoparametric hypersurfaces are homogeneous, translators studied in \cite{F} are given by a function invariant under cohomogenity one actions on the sphere.
On the other hand, Lawn and Ortega \cite{LO} studied the translator given by a function invariant under a cohomegeneity one action on a pseudo-Riemannian manifold.
They showed that the graph of the function gives a translator if and only if the function is a solution of some ODE.  

In this paper, we consider the case where $N$ is a symmetric space $G/K$ of compact type, 
where we give $G/K$ the $G$-invariant metric induced from the $(-1)$-multiple of 
the Killing form of the Lie algebra $\mathfrak g$.  
When the rank of $G/K$ is equal to one (i.e., $G/K=SO(n+1)/SO(n)$, $SU(n+1)/S(U(1)\times U(n))$, $Sp(n+1)/(Sp(1)\times Sp(n))$ or $F_4/{\rm Spin}(9)$), 
we can take the function $r:G/K\to\mathbb{R}$
with $\|\nabla r\|=1$ whose level sets are the orbits of the isotropy group action $K\curvearrowright G/K$.  Then, a function $u$ on a $K$-invariant domain $M$ of $G/K$ which is constant along each orbit of the isotropy group action is given by $u:=V\circ r$ for some function $V$ on $r(M)$.  
It is clear that the shape of the graph of $u$ is dominated by that of $V$.  Hence we suffice to classify the shape of the graph of $V$ to classify that of $u$.  
We obtain the following classification theorem for the shape of the graph of $V$.  

\begin{theorem}\label{thm:vshape}
The graph of $V$ is given by one of the curves obtained by parallel translating curves as in Figures \ref{TypeI}-\ref{TypeV} in the vertical direction.  
The value $\alpha$ in Figures \ref{TypeI}-\ref{TypeV} is the constant given by
$$\alpha=\left\{
\begin{array}{lll}
\displaystyle{\sqrt{\frac{n-1}{2}}\pi} & ({\rm when}\,\,G/K=SO(n+1)/SO(n))\\
\displaystyle{\sqrt{n+1}\pi} & ({\rm when}\,\,G/K=SU(n+1)/S(U(1)\times U(n)))\\
\displaystyle{\sqrt{2(n+2)}\pi} & ({\rm when}\,\,G/K=Sp(n+1)/(Sp(1)\times Sp(n))\\
\displaystyle{\frac{a\pi}{4}} & ({\rm when}\,\,G/K=F_4/{\rm Spin}(9)),
\end{array}
\right.$$
where $a$ is the positive constant such that $a^2$ is equal to the $\frac{1}{4}$-multiple of the maximal sectional curvature of $F_4/{\rm Spin}(9)$.
\begin{figure}[H]
\centering
\scalebox{0.9}{
{\unitlength 0.1in%
\begin{picture}(20.3000,20.8000)(2.0000,-22.0000)%
%
\special{pn 8}%
\special{pa 200 1200}%
\special{pa 2200 1200}%
\special{fp}%
\special{sh 1}%
\special{pa 2200 1200}%
\special{pa 2133 1180}%
\special{pa 2147 1200}%
\special{pa 2133 1220}%
\special{pa 2200 1200}%
\special{fp}%
\put(22.3000,-12.3000){\makebox(0,0)[lb]{$s$}}%
\put(20.3000,-12.8000){\makebox(0,0)[lb]{{\tiny$\alpha$}}}%
\put(3.0000,-13.0000){\makebox(0,0)[lb]{{\tiny $O$}}}%
\special{pn 20}%
\special{pa 632 200}%
\special{pa 635 236}%
\special{pa 640 282}%
\special{pa 645 323}%
\special{pa 650 359}%
\special{pa 655 392}%
\special{pa 660 422}%
\special{pa 665 450}%
\special{pa 670 475}%
\special{pa 675 499}%
\special{pa 680 521}%
\special{pa 685 542}%
\special{pa 690 562}%
\special{pa 700 598}%
\special{pa 705 615}%
\special{pa 710 630}%
\special{pa 715 646}%
\special{pa 725 674}%
\special{pa 735 700}%
\special{pa 745 724}%
\special{pa 760 757}%
\special{pa 770 777}%
\special{pa 775 786}%
\special{pa 780 796}%
\special{pa 785 805}%
\special{pa 790 813}%
\special{pa 795 822}%
\special{pa 810 846}%
\special{pa 815 853}%
\special{pa 820 861}%
\special{pa 830 875}%
\special{pa 835 881}%
\special{pa 840 888}%
\special{pa 845 894}%
\special{pa 850 901}%
\special{pa 865 919}%
\special{pa 870 924}%
\special{pa 875 930}%
\special{pa 880 935}%
\special{pa 885 941}%
\special{pa 905 961}%
\special{pa 910 965}%
\special{pa 920 975}%
\special{pa 935 987}%
\special{pa 940 992}%
\special{pa 950 1000}%
\special{pa 955 1003}%
\special{pa 970 1015}%
\special{pa 980 1021}%
\special{pa 985 1025}%
\special{pa 1000 1034}%
\special{pa 1005 1038}%
\special{pa 1015 1044}%
\special{pa 1020 1046}%
\special{pa 1035 1055}%
\special{pa 1040 1057}%
\special{pa 1045 1060}%
\special{pa 1050 1062}%
\special{pa 1055 1065}%
\special{pa 1065 1069}%
\special{pa 1070 1072}%
\special{pa 1105 1086}%
\special{pa 1110 1087}%
\special{pa 1120 1091}%
\special{pa 1125 1092}%
\special{pa 1135 1096}%
\special{pa 1140 1097}%
\special{pa 1145 1099}%
\special{pa 1155 1101}%
\special{pa 1160 1103}%
\special{pa 1215 1114}%
\special{pa 1220 1114}%
\special{pa 1230 1116}%
\special{pa 1235 1116}%
\special{pa 1240 1117}%
\special{pa 1245 1117}%
\special{pa 1250 1118}%
\special{pa 1255 1118}%
\special{pa 1260 1119}%
\special{pa 1275 1119}%
\special{pa 1280 1120}%
\special{pa 1300 1120}%
\special{fp}%
\special{pn 20}%
\special{pa 1320 1120}%
\special{pa 1325 1119}%
\special{pa 1340 1119}%
\special{pa 1345 1118}%
\special{ip}%
\special{pa 1365 1116}%
\special{pa 1365 1116}%
\special{pa 1370 1116}%
\special{pa 1380 1114}%
\special{pa 1385 1114}%
\special{pa 1390 1113}%
\special{ip}%
\special{pa 1409 1109}%
\special{pa 1434 1104}%
\special{ip}%
\special{pa 1453 1099}%
\special{pa 1455 1099}%
\special{pa 1460 1097}%
\special{pa 1465 1096}%
\special{pa 1475 1092}%
\special{pa 1477 1092}%
\special{ip}%
\special{pa 1496 1086}%
\special{pa 1519 1076}%
\special{ip}%
\special{pa 1538 1068}%
\special{pa 1545 1065}%
\special{pa 1550 1062}%
\special{pa 1555 1060}%
\special{pa 1560 1057}%
\special{ip}%
\special{pa 1578 1047}%
\special{pa 1580 1046}%
\special{pa 1585 1044}%
\special{pa 1595 1038}%
\special{pa 1599 1035}%
\special{ip}%
\special{pa 1616 1024}%
\special{pa 1620 1021}%
\special{pa 1630 1015}%
\special{pa 1636 1010}%
\special{ip}%
\special{pa 1653 998}%
\special{pa 1660 992}%
\special{pa 1665 987}%
\special{pa 1672 982}%
\special{ip}%
\special{pa 1687 968}%
\special{pa 1690 965}%
\special{pa 1695 961}%
\special{pa 1705 951}%
\special{ip}%
\special{pa 1719 937}%
\special{pa 1720 935}%
\special{pa 1725 930}%
\special{pa 1730 924}%
\special{pa 1735 919}%
\special{pa 1736 918}%
\special{ip}%
\special{pa 1748 903}%
\special{pa 1750 901}%
\special{pa 1755 894}%
\special{pa 1760 888}%
\special{pa 1764 883}%
\special{ip}%
\special{pa 1776 867}%
\special{pa 1780 861}%
\special{pa 1785 853}%
\special{pa 1790 846}%
\special{ip}%
\special{pa 1801 829}%
\special{pa 1805 822}%
\special{pa 1810 813}%
\special{pa 1813 808}%
\special{ip}%
\special{pa 1823 790}%
\special{pa 1825 786}%
\special{pa 1830 777}%
\special{pa 1835 768}%
\special{ip}%
\special{pa 1843 750}%
\special{pa 1854 727}%
\special{ip}%
\special{pa 1861 708}%
\special{pa 1865 700}%
\special{pa 1871 685}%
\special{ip}%
\special{pa 1878 666}%
\special{pa 1885 646}%
\special{pa 1886 643}%
\special{ip}%
\special{pa 1892 624}%
\special{pa 1895 615}%
\special{pa 1900 600}%
\special{ip}%
\special{pa 1905 580}%
\special{pa 1910 562}%
\special{pa 1911 556}%
\special{ip}%
\special{pa 1916 537}%
\special{pa 1920 521}%
\special{pa 1922 512}%
\special{ip}%
\special{pa 1926 493}%
\special{pa 1930 475}%
\special{pa 1931 468}%
\special{ip}%
\special{pa 1935 448}%
\special{pa 1940 424}%
\special{ip}%
\special{pa 1943 404}%
\special{pa 1945 392}%
\special{pa 1947 379}%
\special{ip}%
\special{pa 1950 359}%
\special{pa 1950 359}%
\special{pa 1953 335}%
\special{ip}%
\special{pa 1956 315}%
\special{pa 1959 290}%
\special{ip}%
\special{pa 1961 270}%
\special{pa 1964 245}%
\special{ip}%
\special{pa 1966 225}%
\special{pa 1968 200}%
\special{ip}%
\special{pn 20}%
\special{pn 20}%
\special{pa 823 200}%
\special{pa 824 220}%
\special{ip}%
\special{pa 826 245}%
\special{pa 827 265}%
\special{ip}%
\special{pa 829 290}%
\special{pa 830 299}%
\special{pa 831 310}%
\special{ip}%
\special{pa 833 335}%
\special{pa 835 352}%
\special{pa 835 356}%
\special{ip}%
\special{pa 838 381}%
\special{pa 840 398}%
\special{pa 840 401}%
\special{ip}%
\special{pa 843 426}%
\special{pa 845 439}%
\special{pa 846 445}%
\special{ip}%
\special{pa 849 470}%
\special{pa 850 475}%
\special{pa 852 490}%
\special{ip}%
\special{pa 856 515}%
\special{pa 860 535}%
\special{ip}%
\special{pa 864 560}%
\special{pa 865 566}%
\special{pa 868 580}%
\special{ip}%
\special{pa 873 604}%
\special{pa 875 615}%
\special{pa 877 624}%
\special{ip}%
\special{pa 883 648}%
\special{pa 885 657}%
\special{pa 888 668}%
\special{ip}%
\special{pa 894 692}%
\special{pa 895 695}%
\special{pa 900 712}%
\special{ip}%
\special{pa 907 736}%
\special{pa 910 744}%
\special{pa 914 755}%
\special{ip}%
\special{pa 922 778}%
\special{pa 929 797}%
\special{ip}%
\special{pa 939 820}%
\special{pa 940 823}%
\special{pa 945 835}%
\special{pa 947 839}%
\special{ip}%
\special{pa 958 862}%
\special{pa 965 876}%
\special{pa 967 880}%
\special{ip}%
\special{pa 979 902}%
\special{pa 985 912}%
\special{pa 989 919}%
\special{ip}%
\special{pa 1003 940}%
\special{pa 1005 943}%
\special{pa 1015 956}%
\special{ip}%
\special{pa 1030 976}%
\special{pa 1043 992}%
\special{ip}%
\special{pa 1060 1010}%
\special{pa 1074 1024}%
\special{ip}%
\special{pa 1093 1041}%
\special{pa 1109 1053}%
\special{ip}%
\special{pa 1130 1068}%
\special{pa 1147 1078}%
\special{ip}%
\special{pa 1169 1090}%
\special{pa 1175 1092}%
\special{pa 1180 1095}%
\special{pa 1187 1098}%
\special{ip}%
\special{pa 1211 1106}%
\special{pa 1215 1107}%
\special{pa 1220 1109}%
\special{pa 1225 1110}%
\special{pa 1230 1112}%
\special{pa 1230 1112}%
\special{ip}%
\special{pa 1255 1117}%
\special{pa 1260 1117}%
\special{pa 1265 1118}%
\special{pa 1270 1118}%
\special{pa 1275 1119}%
\special{ip}%
\special{pa 1300 1120}%
\special{pa 1315 1120}%
\special{pa 1320 1119}%
\special{pa 1325 1119}%
\special{pa 1330 1118}%
\special{pa 1335 1118}%
\special{pa 1340 1117}%
\special{pa 1345 1117}%
\special{pa 1370 1112}%
\special{pa 1375 1110}%
\special{pa 1380 1109}%
\special{pa 1385 1107}%
\special{pa 1390 1106}%
\special{pa 1395 1104}%
\special{pa 1400 1103}%
\special{pa 1420 1095}%
\special{pa 1425 1092}%
\special{pa 1435 1088}%
\special{pa 1440 1085}%
\special{pa 1445 1083}%
\special{pa 1475 1065}%
\special{pa 1480 1061}%
\special{pa 1485 1058}%
\special{pa 1515 1034}%
\special{pa 1520 1029}%
\special{pa 1525 1025}%
\special{pa 1550 1000}%
\special{pa 1580 964}%
\special{pa 1595 943}%
\special{pa 1600 935}%
\special{pa 1605 928}%
\special{pa 1615 912}%
\special{pa 1635 876}%
\special{pa 1650 846}%
\special{pa 1655 835}%
\special{pa 1660 823}%
\special{pa 1665 812}%
\special{pa 1680 773}%
\special{pa 1685 759}%
\special{pa 1690 744}%
\special{pa 1700 712}%
\special{pa 1705 695}%
\special{pa 1715 657}%
\special{pa 1720 637}%
\special{pa 1725 615}%
\special{pa 1730 591}%
\special{pa 1735 566}%
\special{pa 1740 538}%
\special{pa 1745 508}%
\special{pa 1750 475}%
\special{pa 1755 439}%
\special{pa 1760 398}%
\special{pa 1765 352}%
\special{pa 1770 299}%
\special{pa 1775 236}%
\special{pa 1777 200}%
\special{fp}%
\put(7.0000,-2.5000){\makebox(0,0)[lb]{$V$}}%
%
\special{pn 8}%
\special{pa 400 2200}%
\special{pa 400 200}%
\special{fp}%
\special{sh 1}%
\special{pa 400 200}%
\special{pa 380 267}%
\special{pa 400 253}%
\special{pa 420 267}%
\special{pa 400 200}%
\special{fp}%
%
\special{pn 8}%
\special{pa 2000 2200}%
\special{pa 2000 200}%
\special{dt 0.045}%
\end{picture}}%
}
\caption{Type I}
\label{TypeI}
\end{figure}
\vspace{-0.5cm}
\begin{figure}[H]
\centering
\begin{minipage}{0.45\columnwidth}
\centering
\scalebox{0.9}{
{\unitlength 0.1in%
\begin{picture}(20.3000,20.0000)(2.0000,-22.0000)%
%
\special{pn 8}%
\special{pa 200 1200}%
\special{pa 2200 1200}%
\special{fp}%
\special{sh 1}%
\special{pa 2200 1200}%
\special{pa 2133 1180}%
\special{pa 2147 1200}%
\special{pa 2133 1220}%
\special{pa 2200 1200}%
\special{fp}%
\put(22.3000,-12.3000){\makebox(0,0)[lb]{$s$}}%
\put(3.0000,-13.0000){\makebox(0,0)[lb]{{\tiny $O$}}}%
\special{pn 20}%
\special{pn 20}%
\special{pa 200 494}%
\special{pa 205 472}%
\special{ip}%
\special{pa 211 445}%
\special{pa 215 429}%
\special{pa 216 424}%
\special{ip}%
\special{pa 222 397}%
\special{pa 225 380}%
\special{pa 226 375}%
\special{ip}%
\special{pa 231 348}%
\special{pa 235 327}%
\special{pa 235 326}%
\special{ip}%
\special{pa 240 299}%
\special{pa 240 298}%
\special{pa 243 277}%
\special{ip}%
\special{pa 248 249}%
\special{pa 250 234}%
\special{pa 251 227}%
\special{ip}%
\special{pa 255 200}%
\special{pa 255 200}%
\special{ip}%
\special{pa 565 2200}%
\special{pa 570 2166}%
\special{pa 575 2133}%
\special{pa 580 2102}%
\special{pa 585 2073}%
\special{pa 590 2046}%
\special{pa 595 2020}%
\special{pa 600 1995}%
\special{pa 605 1971}%
\special{pa 610 1948}%
\special{pa 620 1906}%
\special{pa 625 1886}%
\special{pa 630 1867}%
\special{pa 635 1849}%
\special{pa 645 1815}%
\special{pa 655 1783}%
\special{pa 665 1753}%
\special{pa 670 1739}%
\special{pa 690 1687}%
\special{pa 695 1675}%
\special{pa 700 1664}%
\special{pa 705 1652}%
\special{pa 710 1641}%
\special{pa 715 1631}%
\special{pa 720 1620}%
\special{pa 735 1590}%
\special{pa 760 1545}%
\special{pa 765 1537}%
\special{pa 770 1528}%
\special{pa 785 1504}%
\special{pa 790 1497}%
\special{pa 795 1489}%
\special{pa 800 1482}%
\special{pa 805 1474}%
\special{pa 825 1446}%
\special{pa 830 1440}%
\special{pa 835 1433}%
\special{pa 840 1427}%
\special{pa 845 1420}%
\special{pa 855 1408}%
\special{pa 860 1401}%
\special{pa 885 1371}%
\special{pa 890 1366}%
\special{pa 900 1354}%
\special{pa 905 1349}%
\special{pa 910 1343}%
\special{pa 915 1338}%
\special{pa 920 1332}%
\special{pa 925 1327}%
\special{pa 930 1321}%
\special{pa 940 1311}%
\special{pa 945 1305}%
\special{pa 970 1280}%
\special{pa 975 1274}%
\special{pa 1015 1234}%
\special{pa 1020 1230}%
\special{pa 1080 1170}%
\special{pa 1085 1166}%
\special{pa 1125 1126}%
\special{pa 1130 1120}%
\special{pa 1155 1095}%
\special{pa 1160 1089}%
\special{pa 1170 1079}%
\special{pa 1175 1073}%
\special{pa 1180 1068}%
\special{pa 1185 1062}%
\special{pa 1190 1057}%
\special{pa 1195 1051}%
\special{pa 1200 1046}%
\special{pa 1210 1034}%
\special{pa 1215 1029}%
\special{pa 1240 999}%
\special{pa 1245 992}%
\special{pa 1255 980}%
\special{pa 1260 973}%
\special{pa 1265 967}%
\special{pa 1270 960}%
\special{pa 1275 954}%
\special{pa 1295 926}%
\special{pa 1300 918}%
\special{pa 1305 911}%
\special{pa 1310 903}%
\special{pa 1315 896}%
\special{pa 1330 872}%
\special{pa 1335 863}%
\special{pa 1340 855}%
\special{pa 1365 810}%
\special{pa 1380 780}%
\special{pa 1385 769}%
\special{pa 1390 759}%
\special{pa 1395 748}%
\special{pa 1400 736}%
\special{pa 1405 725}%
\special{pa 1410 713}%
\special{pa 1430 661}%
\special{pa 1435 647}%
\special{pa 1445 617}%
\special{pa 1455 585}%
\special{pa 1465 551}%
\special{pa 1470 533}%
\special{pa 1475 514}%
\special{pa 1480 494}%
\special{pa 1490 452}%
\special{pa 1495 429}%
\special{pa 1500 405}%
\special{pa 1505 380}%
\special{pa 1510 354}%
\special{pa 1515 327}%
\special{pa 1520 298}%
\special{pa 1525 267}%
\special{pa 1530 234}%
\special{pa 1535 200}%
\special{fp}%
\special{pn 20}%
\special{pa 1848 2180}%
\special{pa 1850 2166}%
\special{pa 1852 2155}%
\special{ip}%
\special{pa 1855 2135}%
\special{pa 1855 2133}%
\special{pa 1859 2110}%
\special{ip}%
\special{pa 1862 2090}%
\special{pa 1865 2073}%
\special{pa 1866 2065}%
\special{ip}%
\special{pa 1870 2045}%
\special{pa 1875 2021}%
\special{ip}%
\special{pa 1879 2001}%
\special{pa 1880 1995}%
\special{pa 1884 1976}%
\special{ip}%
\special{pa 1888 1956}%
\special{pa 1890 1948}%
\special{pa 1894 1932}%
\special{ip}%
\special{pa 1899 1912}%
\special{pa 1900 1906}%
\special{pa 1905 1888}%
\special{ip}%
\special{pa 1910 1868}%
\special{pa 1910 1867}%
\special{pa 1915 1849}%
\special{pa 1917 1844}%
\special{ip}%
\special{pa 1922 1824}%
\special{pa 1925 1815}%
\special{pa 1930 1800}%
\special{ip}%
\special{pa 1936 1781}%
\special{pa 1944 1757}%
\special{ip}%
\special{pa 1950 1738}%
\special{pa 1960 1714}%
\special{ip}%
\special{pa 1967 1695}%
\special{pa 1970 1687}%
\special{pa 1975 1675}%
\special{pa 1976 1672}%
\special{ip}%
\special{pa 1984 1653}%
\special{pa 1985 1652}%
\special{pa 1990 1641}%
\special{pa 1995 1631}%
\special{pa 1995 1631}%
\special{ip}%
\special{pa 2004 1612}%
\special{pa 2015 1590}%
\special{pa 2015 1590}%
\special{ip}%
\special{pa 2025 1572}%
\special{pa 2037 1550}%
\special{ip}%
\special{pa 2047 1533}%
\special{pa 2050 1528}%
\special{pa 2061 1511}%
\special{ip}%
\special{pa 2072 1494}%
\special{pa 2075 1489}%
\special{pa 2080 1482}%
\special{pa 2085 1474}%
\special{pa 2086 1473}%
\special{ip}%
\special{pa 2097 1457}%
\special{pa 2105 1446}%
\special{pa 2110 1440}%
\special{pa 2113 1436}%
\special{ip}%
\special{pa 2125 1420}%
\special{pa 2125 1420}%
\special{pa 2135 1408}%
\special{pa 2140 1401}%
\special{pa 2140 1401}%
\special{ip}%
\special{pa 2153 1385}%
\special{pa 2165 1371}%
\special{pa 2170 1366}%
\special{ip}%
\special{pa 2183 1351}%
\special{pa 2185 1349}%
\special{pa 2190 1343}%
\special{pa 2195 1338}%
\special{pa 2200 1332}%
\special{ip}%
\put(13.4000,-4.0000){\makebox(0,0)[lb]{$V$}}%
%
\special{pn 8}%
\special{pa 400 2200}%
\special{pa 400 200}%
\special{fp}%
\special{sh 1}%
\special{pa 400 200}%
\special{pa 380 267}%
\special{pa 400 253}%
\special{pa 420 267}%
\special{pa 400 200}%
\special{fp}%
%
\special{pn 8}%
\special{pa 2000 2200}%
\special{pa 2000 200}%
\special{dt 0.045}%
\put(20.3000,-12.8000){\makebox(0,0)[lb]{{\tiny$\alpha$}}}%
\end{picture}}%
}
\caption{Type II}
\end{minipage}
\begin{minipage}{0.45\columnwidth}
\centering
\scalebox{0.9}{
{\unitlength 0.1in%
\begin{picture}(20.5000,20.0000)(1.8000,-22.0000)%
%
\special{pn 8}%
\special{pa 200 1200}%
\special{pa 2200 1200}%
\special{fp}%
\special{sh 1}%
\special{pa 2200 1200}%
\special{pa 2133 1180}%
\special{pa 2147 1200}%
\special{pa 2133 1220}%
\special{pa 2200 1200}%
\special{fp}%
\put(22.3000,-12.3000){\makebox(0,0)[lb]{$s$}}%
\put(20.3000,-12.8000){\makebox(0,0)[lb]{{\tiny$\alpha$}}}%
\put(1.8000,-11.7000){\makebox(0,0)[lb]{{\tiny $O$}}}%
\special{pn 20}%
\special{pn 20}%
\special{pa 200 1232}%
\special{pa 205 1238}%
\special{pa 210 1243}%
\special{pa 214 1247}%
\special{ip}%
\special{pa 230 1266}%
\special{pa 235 1271}%
\special{pa 244 1282}%
\special{ip}%
\special{pa 260 1301}%
\special{pa 265 1308}%
\special{pa 273 1317}%
\special{ip}%
\special{pa 288 1338}%
\special{pa 290 1340}%
\special{pa 295 1346}%
\special{pa 301 1354}%
\special{ip}%
\special{pa 315 1375}%
\special{pa 320 1382}%
\special{pa 325 1389}%
\special{pa 327 1392}%
\special{ip}%
\special{pa 341 1413}%
\special{pa 350 1428}%
\special{pa 351 1430}%
\special{ip}%
\special{pa 364 1452}%
\special{pa 374 1470}%
\special{ip}%
\special{pa 386 1492}%
\special{pa 395 1511}%
\special{ip}%
\special{pa 406 1534}%
\special{pa 410 1541}%
\special{pa 415 1552}%
\special{pa 415 1552}%
\special{ip}%
\special{pa 425 1576}%
\special{pa 430 1587}%
\special{pa 433 1594}%
\special{ip}%
\special{pa 442 1618}%
\special{pa 449 1637}%
\special{ip}%
\special{pa 458 1661}%
\special{pa 464 1681}%
\special{ip}%
\special{pa 472 1705}%
\special{pa 475 1715}%
\special{pa 478 1725}%
\special{ip}%
\special{pa 485 1749}%
\special{pa 490 1767}%
\special{pa 490 1769}%
\special{ip}%
\special{pa 497 1793}%
\special{pa 500 1806}%
\special{pa 502 1813}%
\special{ip}%
\special{pa 508 1838}%
\special{pa 510 1848}%
\special{pa 512 1858}%
\special{ip}%
\special{pa 517 1883}%
\special{pa 520 1895}%
\special{pa 522 1903}%
\special{ip}%
\special{pa 527 1928}%
\special{pa 530 1946}%
\special{pa 530 1948}%
\special{ip}%
\special{pa 535 1973}%
\special{pa 535 1973}%
\special{pa 538 1993}%
\special{ip}%
\special{pa 543 2018}%
\special{pa 545 2033}%
\special{pa 546 2038}%
\special{ip}%
\special{pa 550 2064}%
\special{pa 550 2066}%
\special{pa 553 2084}%
\special{ip}%
\special{pa 556 2109}%
\special{pa 559 2129}%
\special{ip}%
\special{pa 562 2154}%
\special{pa 565 2175}%
\special{ip}%
\special{ip}%
\special{pa 880 200}%
\special{pa 890 254}%
\special{pa 895 280}%
\special{pa 900 305}%
\special{pa 905 329}%
\special{pa 910 352}%
\special{pa 920 394}%
\special{pa 925 414}%
\special{pa 930 433}%
\special{pa 935 451}%
\special{pa 945 485}%
\special{pa 955 517}%
\special{pa 965 547}%
\special{pa 970 561}%
\special{pa 990 613}%
\special{pa 995 625}%
\special{pa 1000 636}%
\special{pa 1005 648}%
\special{pa 1010 659}%
\special{pa 1015 669}%
\special{pa 1020 680}%
\special{pa 1035 710}%
\special{pa 1060 755}%
\special{pa 1065 763}%
\special{pa 1070 772}%
\special{pa 1085 796}%
\special{pa 1090 803}%
\special{pa 1095 811}%
\special{pa 1100 818}%
\special{pa 1105 826}%
\special{pa 1125 854}%
\special{pa 1130 860}%
\special{pa 1135 867}%
\special{pa 1140 873}%
\special{pa 1145 880}%
\special{pa 1155 892}%
\special{pa 1160 899}%
\special{pa 1185 929}%
\special{pa 1190 934}%
\special{pa 1200 946}%
\special{pa 1205 951}%
\special{pa 1210 957}%
\special{pa 1215 962}%
\special{pa 1220 968}%
\special{pa 1225 973}%
\special{pa 1230 979}%
\special{pa 1240 989}%
\special{pa 1245 995}%
\special{pa 1270 1020}%
\special{pa 1275 1026}%
\special{pa 1315 1066}%
\special{pa 1320 1070}%
\special{pa 1380 1130}%
\special{pa 1385 1134}%
\special{pa 1425 1174}%
\special{pa 1430 1180}%
\special{pa 1455 1205}%
\special{pa 1460 1211}%
\special{pa 1470 1221}%
\special{pa 1475 1227}%
\special{pa 1480 1232}%
\special{pa 1485 1238}%
\special{pa 1490 1243}%
\special{pa 1495 1249}%
\special{pa 1500 1254}%
\special{pa 1510 1266}%
\special{pa 1515 1271}%
\special{pa 1540 1301}%
\special{pa 1545 1308}%
\special{pa 1555 1320}%
\special{pa 1560 1327}%
\special{pa 1565 1333}%
\special{pa 1570 1340}%
\special{pa 1575 1346}%
\special{pa 1595 1374}%
\special{pa 1600 1382}%
\special{pa 1605 1389}%
\special{pa 1610 1397}%
\special{pa 1615 1404}%
\special{pa 1630 1428}%
\special{pa 1635 1437}%
\special{pa 1640 1445}%
\special{pa 1665 1490}%
\special{pa 1680 1520}%
\special{pa 1685 1531}%
\special{pa 1690 1541}%
\special{pa 1695 1552}%
\special{pa 1700 1564}%
\special{pa 1705 1575}%
\special{pa 1710 1587}%
\special{pa 1730 1639}%
\special{pa 1735 1653}%
\special{pa 1745 1683}%
\special{pa 1755 1715}%
\special{pa 1765 1749}%
\special{pa 1770 1767}%
\special{pa 1775 1786}%
\special{pa 1780 1806}%
\special{pa 1790 1848}%
\special{pa 1795 1871}%
\special{pa 1800 1895}%
\special{pa 1805 1920}%
\special{pa 1810 1946}%
\special{pa 1815 1973}%
\special{pa 1820 2002}%
\special{pa 1825 2033}%
\special{pa 1830 2066}%
\special{pa 1835 2100}%
\special{pa 1840 2137}%
\special{pa 1845 2176}%
\special{pa 1848 2200}%
\special{fp}%
\special{pn 20}%
\special{pa 2164 222}%
\special{pa 2169 249}%
\special{ip}%
\special{pa 2173 270}%
\special{pa 2175 280}%
\special{pa 2178 297}%
\special{ip}%
\special{pa 2183 319}%
\special{pa 2185 329}%
\special{pa 2189 346}%
\special{ip}%
\special{pa 2194 367}%
\special{pa 2200 394}%
\special{ip}%
\put(9.7000,-4.0000){\makebox(0,0)[lb]{$V$}}%
%
\special{pn 8}%
\special{pa 400 2200}%
\special{pa 400 200}%
\special{fp}%
\special{sh 1}%
\special{pa 400 200}%
\special{pa 380 267}%
\special{pa 400 253}%
\special{pa 420 267}%
\special{pa 400 200}%
\special{fp}%
%
\special{pn 8}%
\special{pa 2000 2200}%
\special{pa 2000 200}%
\special{dt 0.045}%
\end{picture}}%
}
\caption{Type III}
\end{minipage}
\end{figure}
\begin{figure}[H]
\centering
\begin{minipage}{0.45\columnwidth}
\centering
\scalebox{0.9}{
{\unitlength 0.1in%
\begin{picture}(20.3000,20.0000)(2.0000,-22.0000)%
%
\special{pn 8}%
\special{pa 200 1200}%
\special{pa 2200 1200}%
\special{fp}%
\special{sh 1}%
\special{pa 2200 1200}%
\special{pa 2133 1180}%
\special{pa 2147 1200}%
\special{pa 2133 1220}%
\special{pa 2200 1200}%
\special{fp}%
\put(22.3000,-12.3000){\makebox(0,0)[lb]{$s$}}%
\put(20.3000,-12.8000){\makebox(0,0)[lb]{{\tiny$\alpha$}}}%
\put(3.0000,-13.0000){\makebox(0,0)[lb]{{\tiny$O$}}}%
\special{pn 20}%
\special{pn 20}%
\special{pa 200 1154}%
\special{pa 205 1154}%
\special{pa 210 1153}%
\special{pa 220 1153}%
\special{pa 222 1153}%
\special{ip}%
\special{pa 250 1151}%
\special{pa 255 1150}%
\special{pa 260 1150}%
\special{pa 265 1149}%
\special{pa 272 1149}%
\special{ip}%
\special{pa 300 1146}%
\special{pa 310 1146}%
\special{pa 315 1145}%
\special{pa 320 1145}%
\special{pa 322 1145}%
\special{ip}%
\special{pa 350 1142}%
\special{pa 360 1140}%
\special{pa 365 1140}%
\special{pa 370 1139}%
\special{pa 372 1139}%
\special{ip}%
\special{ip}%
\special{pa 400 1136}%
\special{pa 410 1134}%
\special{pa 415 1134}%
\special{pa 425 1132}%
\special{pa 430 1132}%
\special{pa 440 1130}%
\special{pa 445 1130}%
\special{pa 455 1128}%
\special{pa 460 1128}%
\special{pa 475 1125}%
\special{pa 480 1125}%
\special{pa 495 1122}%
\special{pa 500 1122}%
\special{pa 525 1117}%
\special{pa 530 1117}%
\special{pa 575 1108}%
\special{pa 580 1108}%
\special{pa 625 1099}%
\special{pa 630 1097}%
\special{pa 675 1088}%
\special{pa 680 1086}%
\special{pa 700 1082}%
\special{pa 705 1080}%
\special{pa 720 1077}%
\special{pa 725 1075}%
\special{pa 740 1072}%
\special{pa 745 1070}%
\special{pa 755 1068}%
\special{pa 760 1066}%
\special{pa 765 1065}%
\special{pa 770 1063}%
\special{pa 780 1061}%
\special{pa 785 1059}%
\special{pa 790 1058}%
\special{pa 795 1056}%
\special{pa 800 1055}%
\special{pa 805 1053}%
\special{pa 810 1052}%
\special{pa 815 1050}%
\special{pa 820 1049}%
\special{pa 825 1047}%
\special{pa 830 1046}%
\special{pa 835 1044}%
\special{pa 840 1043}%
\special{pa 845 1041}%
\special{pa 850 1040}%
\special{pa 860 1036}%
\special{pa 865 1035}%
\special{pa 875 1031}%
\special{pa 880 1030}%
\special{pa 890 1026}%
\special{pa 895 1025}%
\special{pa 915 1017}%
\special{pa 920 1016}%
\special{pa 1015 978}%
\special{pa 1020 975}%
\special{pa 1035 969}%
\special{pa 1040 966}%
\special{pa 1050 962}%
\special{pa 1055 959}%
\special{pa 1065 955}%
\special{pa 1070 952}%
\special{pa 1075 950}%
\special{pa 1080 947}%
\special{pa 1085 945}%
\special{pa 1090 942}%
\special{pa 1095 940}%
\special{pa 1105 934}%
\special{pa 1110 932}%
\special{pa 1120 926}%
\special{pa 1125 924}%
\special{pa 1140 915}%
\special{pa 1145 913}%
\special{pa 1190 886}%
\special{pa 1195 882}%
\special{pa 1210 873}%
\special{pa 1215 869}%
\special{pa 1225 863}%
\special{pa 1230 859}%
\special{pa 1235 856}%
\special{pa 1240 852}%
\special{pa 1245 849}%
\special{pa 1255 841}%
\special{pa 1260 838}%
\special{pa 1315 794}%
\special{pa 1320 789}%
\special{pa 1325 785}%
\special{pa 1330 780}%
\special{pa 1335 776}%
\special{pa 1345 766}%
\special{pa 1350 762}%
\special{pa 1380 732}%
\special{pa 1385 726}%
\special{pa 1390 721}%
\special{pa 1395 715}%
\special{pa 1400 710}%
\special{pa 1435 668}%
\special{pa 1440 661}%
\special{pa 1445 655}%
\special{pa 1465 627}%
\special{pa 1470 619}%
\special{pa 1475 612}%
\special{pa 1495 580}%
\special{pa 1500 571}%
\special{pa 1505 563}%
\special{pa 1515 545}%
\special{pa 1540 495}%
\special{pa 1545 484}%
\special{pa 1550 472}%
\special{pa 1555 461}%
\special{pa 1560 448}%
\special{pa 1565 436}%
\special{pa 1570 423}%
\special{pa 1580 395}%
\special{pa 1590 365}%
\special{pa 1595 349}%
\special{pa 1600 332}%
\special{pa 1610 296}%
\special{pa 1615 276}%
\special{pa 1620 255}%
\special{pa 1625 232}%
\special{pa 1630 208}%
\special{pa 1632 200}%
\special{fp}%
\put(14.0000,-4.0000){\makebox(0,0)[lb]{$V$}}%
%
\special{pn 8}%
\special{pa 400 2200}%
\special{pa 400 200}%
\special{fp}%
\special{sh 1}%
\special{pa 400 200}%
\special{pa 380 267}%
\special{pa 400 253}%
\special{pa 420 267}%
\special{pa 400 200}%
\special{fp}%
%
\special{pn 8}%
\special{pa 2000 2200}%
\special{pa 2000 200}%
\special{dt 0.045}%
\end{picture}}%
}
\caption{Type IV}
\end{minipage}
\begin{minipage}{0.45\columnwidth}
\centering
\scalebox{0.9}{
{\unitlength 0.1in%
\begin{picture}(20.3000,20.0000)(2.0000,-22.0000)%
%
\special{pn 8}%
\special{pa 200 1200}%
\special{pa 2200 1200}%
\special{fp}%
\special{sh 1}%
\special{pa 2200 1200}%
\special{pa 2133 1180}%
\special{pa 2147 1200}%
\special{pa 2133 1220}%
\special{pa 2200 1200}%
\special{fp}%
\put(22.3000,-12.3000){\makebox(0,0)[lb]{$s$}}%
\put(20.3000,-12.9000){\makebox(0,0)[lb]{{\tiny$\alpha$}}}%
\put(3.0000,-13.0000){\makebox(0,0)[lb]{{\tiny$O$}}}%
\special{pn 20}%
\special{pa 768 200}%
\special{pa 770 208}%
\special{pa 775 232}%
\special{pa 780 255}%
\special{pa 785 276}%
\special{pa 790 296}%
\special{pa 800 332}%
\special{pa 805 349}%
\special{pa 810 365}%
\special{pa 820 395}%
\special{pa 830 423}%
\special{pa 835 436}%
\special{pa 840 448}%
\special{pa 845 461}%
\special{pa 850 472}%
\special{pa 855 484}%
\special{pa 860 495}%
\special{pa 885 545}%
\special{pa 895 563}%
\special{pa 900 571}%
\special{pa 905 580}%
\special{pa 925 612}%
\special{pa 930 619}%
\special{pa 935 627}%
\special{pa 955 655}%
\special{pa 960 661}%
\special{pa 965 668}%
\special{pa 1000 710}%
\special{pa 1005 715}%
\special{pa 1010 721}%
\special{pa 1015 726}%
\special{pa 1020 732}%
\special{pa 1050 762}%
\special{pa 1055 766}%
\special{pa 1065 776}%
\special{pa 1070 780}%
\special{pa 1075 785}%
\special{pa 1080 789}%
\special{pa 1085 794}%
\special{pa 1140 838}%
\special{pa 1145 841}%
\special{pa 1155 849}%
\special{pa 1160 852}%
\special{pa 1165 856}%
\special{pa 1170 859}%
\special{pa 1175 863}%
\special{pa 1185 869}%
\special{pa 1190 873}%
\special{pa 1205 882}%
\special{pa 1210 886}%
\special{pa 1255 913}%
\special{pa 1260 915}%
\special{pa 1275 924}%
\special{pa 1280 926}%
\special{pa 1290 932}%
\special{pa 1295 934}%
\special{pa 1305 940}%
\special{pa 1310 942}%
\special{pa 1315 945}%
\special{pa 1320 947}%
\special{pa 1325 950}%
\special{pa 1330 952}%
\special{pa 1335 955}%
\special{pa 1345 959}%
\special{pa 1350 962}%
\special{pa 1360 966}%
\special{pa 1365 969}%
\special{pa 1380 975}%
\special{pa 1385 978}%
\special{pa 1480 1016}%
\special{pa 1485 1017}%
\special{pa 1505 1025}%
\special{pa 1510 1026}%
\special{pa 1520 1030}%
\special{pa 1525 1031}%
\special{pa 1535 1035}%
\special{pa 1540 1036}%
\special{pa 1550 1040}%
\special{pa 1555 1041}%
\special{pa 1560 1043}%
\special{pa 1565 1044}%
\special{pa 1570 1046}%
\special{pa 1575 1047}%
\special{pa 1580 1049}%
\special{pa 1585 1050}%
\special{pa 1590 1052}%
\special{pa 1595 1053}%
\special{pa 1600 1055}%
\special{pa 1605 1056}%
\special{pa 1610 1058}%
\special{pa 1615 1059}%
\special{pa 1620 1061}%
\special{pa 1630 1063}%
\special{pa 1635 1065}%
\special{pa 1640 1066}%
\special{pa 1645 1068}%
\special{pa 1655 1070}%
\special{pa 1660 1072}%
\special{pa 1675 1075}%
\special{pa 1680 1077}%
\special{pa 1695 1080}%
\special{pa 1700 1082}%
\special{pa 1720 1086}%
\special{pa 1725 1088}%
\special{pa 1770 1097}%
\special{pa 1775 1099}%
\special{pa 1820 1108}%
\special{pa 1825 1108}%
\special{pa 1870 1117}%
\special{pa 1875 1117}%
\special{pa 1900 1122}%
\special{pa 1905 1122}%
\special{pa 1920 1125}%
\special{pa 1925 1125}%
\special{pa 1940 1128}%
\special{pa 1945 1128}%
\special{pa 1955 1130}%
\special{pa 1960 1130}%
\special{pa 1970 1132}%
\special{pa 1975 1132}%
\special{pa 1985 1134}%
\special{pa 1990 1134}%
\special{pa 2000 1136}%
\special{fp}%
\special{pn 20}%
\special{pa 2022 1138}%
\special{pa 2025 1139}%
\special{pa 2030 1139}%
\special{pa 2035 1140}%
\special{pa 2040 1140}%
\special{pa 2050 1142}%
\special{ip}%
\special{pa 2072 1144}%
\special{pa 2075 1144}%
\special{pa 2080 1145}%
\special{pa 2085 1145}%
\special{pa 2090 1146}%
\special{pa 2100 1146}%
\special{ip}%
\special{pa 2122 1148}%
\special{pa 2125 1149}%
\special{pa 2135 1149}%
\special{pa 2140 1150}%
\special{pa 2145 1150}%
\special{pa 2150 1151}%
\special{ip}%
\special{pa 2172 1152}%
\special{pa 2175 1152}%
\special{pa 2180 1153}%
\special{pa 2190 1153}%
\special{pa 2195 1154}%
\special{pa 2200 1154}%
\special{ip}%
\put(8.5000,-3.8000){\makebox(0,0)[lb]{$V$}}%
%
\special{pn 8}%
\special{pa 400 2200}%
\special{pa 400 200}%
\special{fp}%
\special{sh 1}%
\special{pa 400 200}%
\special{pa 380 267}%
\special{pa 400 253}%
\special{pa 420 267}%
\special{pa 400 200}%
\special{fp}%
%
\special{pn 8}%
\special{pa 2000 2200}%
\special{pa 2000 200}%
\special{dt 0.045}%
\end{picture}}%
}
\caption{Type V}
\label{TypeV}
\end{minipage}
\end{figure}
\end{theorem}

Next we consider the case where $G/K$ is a higher rank irreducible symmetric space of compact type and a translator (for the mean curvature flow) given by 
a graph of a function on $G/K$ which is invariant under a Hermann action $H\curvearrowright G/K$ of cohomegeneity two, where 
{\it Hermann action} means that $H$ is a symmetric subgroup of $G$.  Assume that $H\curvearrowright G/K$ is commutative, 
that is, $\theta_K\circ\theta_H=\theta_H\circ\theta_K$ holds for the involutions $\theta_K$ and $\theta_H$ of $G$ satisfying 
$({\rm Fix}\,\theta_K)_0\subset K\subset{\rm Fix}\,\theta_K$ and $({\rm Fix}\,\theta_H)_0\subset H\subset{\rm Fix}\,\theta_H$, 
where ${\rm Fix}(\cdot)$ is the fixed point group of $(\cdot)$ and $(\cdot)_0$ is the identity component of $(\cdot)$.  
Here we note that Hermann actions are hyperpolar actions, where a {\it hyperpolar action} means an isometric action of a compact Lie group on $G/K$ which admits 
a complete flat totally geodesic submanifold in $G/K$ meeting all orbits of the action orthogonally.  
The complete flat totally geodesic submanifold is called a {\it flat section} of this action.

Let $r=(r_1,r_2):G/K\to\mathbb{R}^2$ be a map on $G/K$ with $g(\nabla r_i,\nabla r_j)=\delta_{ij}$ ($i,j\in\{1,2\}$) whose level sets are the orbits of 
the action $H\curvearrowright G/K$.  
Then, a function $u$ on a $H$-invariant domain $M$ of $G/K$ is invariant under the action $H\curvearrowright G/K$ if and only if $u$ is described as 
$u=V\circ r$ for some function $V$ on $r(M)$.  Let $\Sigma$ be the flat section of the $H$-action through $o:=eK$, were we note that $\Sigma$ is 
diffeomorphic to a torus $T^2(=S^1\times S^1)$.  
Let $\mathfrak g=\mathfrak k\oplus\mathfrak p$ be the canonical decomposition associated to the symmetric pair $(G,K)$.  The space $\mathfrak p$ is identified with 
the tangent space $T_o(G/K)$ through the restriction $\pi_{\ast e}|_{\mathfrak p}$ of the differential $\pi_{\ast e}$ of the natural projection $\pi:G\to G/K$.  
There exists the maximal abelian subspace $\mathfrak a$ of $\mathfrak p$ satisfying $\exp_o(\mathfrak a)=\Sigma$, where $\exp_o$ is the exponential map of $G/K$ at $o$.  
Let $\mathcal C(\subset\mathfrak a)$ be a Weyl domain and $W$ be the Weyl group.
Denote by $\vX$ the tangent vector field on $\exp_o(\mathcal C)$ defined by assigning the mean curvature vector of the orbit $H\cdot w$ at $w$ to each $w\in\exp_o(\mathcal C)$.
By identifying $\exp_o(\mathcal C)$ with $\mathcal C$, we regard $\vX$ as a tangent vector field on $\mathcal C$.  

\begin{theorem}\label{thm:pde-cohomo2}
The graph $\Gamma$ of $u=V\circ r$ is a translator if and only if $V$ satisfies 
\begin{equation}
\sum_{i,j=1}^2\frac{\partial^2V}{\partial x_i\partial x_j}\frac{\partial V}{\partial x_i}\frac{\partial V}{\partial x_j}
-\left(1+|\nabla V|^2\right)\left(\sum_{i=1}^2X_i\frac{\partial V}{\partial x_i}+\Delta V-1\right)=0,\label{eq:pde1-cohomo2}
\end{equation}
where $(x_1,x_2)$ is the Euclidean coordinate of $\mathfrak a$, $V$ is regarded as a function on $\mathfrak a$ through $(x_1,x_2):\mathfrak a\to\mathbb R^2$ 
and $X_i$ ($i=1,2$) are the components of the tangent vector field $\vX$ on $\mathcal C$ with respect to the Euclidean coordinate $(x_1,x_2)$ of $\mathfrak a$ 
(i.e., $\vX=\sum\limits_{i=1}^2X_i\frac{\partial}{\partial x_i}$).  
\end{theorem}

For all commutative Hermann actions of cohomogeneity two on an irreducible symmetric space of compact type, the explicit descriptions of the component 
$(X_1,X_2)$ of the above tangent vector field $\vX$ on $\mathcal C$ are given in \cite{K}.  By using the explicit description, we can describe the PDE 
{\rm (\ref{eq:pde1-cohomo2})} of order two explicitly.  
Clearly we can choose the above function $r$ as $r|_{\exp_o(\mathcal C)}=(x_1,x_2)\circ(\exp_o|_{\mathcal C})^{-1}$ holds.  
According to \cite{K}, $\vX$ is described as $\vX=\nabla\rho$ for some convex function $\rho$ on $\mathcal C$.  
Next we consider the case where $V$ is constant along each level set of $\rho$.  
In this case, $\nabla V$ is described as $\nabla V=F\vX$ for some function $F$ on $\mathcal C$.  
It is clear that the shape of the graph of $u$ is dominated by $F$.  Hence we suffice to investigate $F$ to classify the shape of the graph of $u$.
In this case, we obtain the following fact.  

\begin{theorem}\label{thm:ode-cohomo2}
Assume that $V$ is constant along each level set of $\rho$ and let $F$ be the function on $\mathcal C$ satisfying $\nabla V=F\vX$.  
Then the graph $\Gamma$ of $u=V\circ r$ is a translator if and only if $F$ satisfies 
\begin{equation}\label{eq:pde2-cohomo2}
\vX(F)=\frac{1}{2}\vX(|\vX|^2)F^3-(1+|\vX|^2F^2)((|\vX|^2+{\rm div}\,\vX)F-1).
\end{equation}
\end{theorem}

By using the explicit descriptions of $\vX$ in \cite{K}, we can describe the PDE {\rm (\ref{eq:pde2-cohomo2})} of order one explicitly.  
As one cexample we investigate the shape of the graph of $V$ in the case where the Hermann action $H\curvearrowright G/K$ is the dual action of the Hermann type action $SO_0(1,2)\curvearrowright SL(3,\mathbb R)/SO(3)$.

In Section 2, we investigate translators which are invariant under the isotropy group action of rank one symmetric spaces of compact type and prove Theorem 1.1.  
In Section 3, we investigate translators which are invariant under Hermann action of cohomogeneity two on higher rank symmetric spaces of compact type and prove Theorems 1.2 
and 1.3.  

\section{The case of cohomogenity one} 
Let $(N,g)$ be an $n$-dimensional Riemannian manifold and $u:M\to\mathbb{R}$ be a ($C^{\infty}$)-function on a domain $M$ of $N$.  
Let $f$ be the graph embedding of $u$, that is, the embedding $f:M\hookrightarrow(N,g)$ defined by $f(x)=(x,u(x))\,\,\,(x\in M)$.  
Denote by $\Gamma$ the graph of $u$ and $H$ the mean curvature vector field of $f$.  
Also, denote by $\nabla(\cdot)$ and ${\rm div}(\cdot)$ be the gradient vector field and the divergence of $(\cdot)$ with respect to $g$, respectively.  
For the translatority of $\Gamma$, the following fact holds (see \cite{CSS}, \cite{F} and \cite{LO}).  

\begin{lemma}\label{lemma:graph-soliton-condition}
If $\Gamma$ is a translator, $u$ satisfies
\begin{equation}
\sqrt{1+\|\nabla u\|^2}~{\rm div}\left(\frac{\nabla u}{\sqrt{1+\|\nabla u\|^2}}\right)=1.\label{eq:graph-soliton}
\end{equation}
Conversely, if $u$ satisfies {\rm (\ref{eq:graph-soliton})}, $\Gamma$ is a translator.
\end{lemma}

We consider the case where $(N,g)$ is a rank one symmetric space $G/K$ of compact type and $u$ is a function on a $K$-invariant domain $M$ of $G/K$ which is invariant under the isotropy group action $K\curvearrowright G/K$, 
where we give $G/K$ the $G$-invariant metric induced from the $(-1)$-multiple of the Killing form of $\mathfrak g$.  
Let $r:G/K\to\mathbb{R}$ be the function on $G/K$ with $\|\nabla r\|=1$ whose level sets are the orbits of $K\curvearrowright G/K$.  Then, since $u$ is invariant under 
the action $K\curvearrowright G/K$, $u$ is described as $u=V\circ r$ for some function $V$ on $r(M)$.  
According to the result by Lawn and Ortega (Theorem 3.5 of \cite{LO}) for the graph of a function on a pseudo-Riemaniann manifold which is invariant under 
a cohomogeneity one proper isometric action of a Lie group, we obtain the following fact.  

\begin{proposition}\label{prop:ode-cohomo1}
The graph $\Gamma$ of $V\circ r$ is a translator if and only if $V$ satisfies 
\begin{equation}\label{eq:ode-cohomo1}
V''(s)=\left(1+V'(s)^2\right)\left(1-h(s)V'(s)\right)
\end{equation}
where $(\cdot)'$ denotes the derivative of $(\cdot)$ and $h(s)$ is the constant mean curvature of the orbit $r^{-1}(s)$ of $K\curvearrowright G/K$.
\end{proposition}

Let $\mathfrak p$ be as in Introduction. $\mathfrak a$ be a maximal abelian subspace of $\mathfrak p$ and $\triangle_+$ be the positive root system with respect to 
$\mathfrak a$.  Then we have 
\begin{equation}\label{eq:root-system}
\triangle_+=\left\{\begin{array}{ll}
\displaystyle{\left\{\frac{1}{\sqrt{2(n-1)}}\,\langle\ve,\cdot\rangle\right\}} & ({\rm when}\,\,G/K=SO(n+1)/SO(n))\\
\displaystyle{\left\{\frac{1}{4\sqrt{n+1}}\,\langle\ve,\cdot\rangle,\,\frac{1}{2\sqrt{n+1}}\,\langle\ve,\cdot\rangle\right\}} & ({\rm when}\,\,G/K=SU(n+1)/S(U(1)\times U(n)))\\
\displaystyle{\left\{\frac{1}{4\sqrt{2(n+2)}}\,\langle\ve,\cdot\rangle,\,\frac{1}{2\sqrt{2(n+2)}}\,\langle\ve,\cdot\rangle\right\}} & 
({\rm when}\,\,G/K=Sp(n+1)/(Sp(1)\times Sp(n)))\\
\displaystyle{\{a\,\langle\ve,\cdot\rangle,\,2a\,\langle\ve,\cdot\rangle\}} & ({\rm when}\,\,G/K=F_4/{\rm Spin}(9)),
\end{array}\right.
\end{equation}
where $\ve$ is a unit normal vector of $\mathfrak a$, $\langle\,\,,\,\,\rangle$ is the restriction of the $(-1)$-multiple of of the Killing form to $\mathfrak a$ 
and $a$ is the positive constant stated in Theorem 1.1.  
Thus, in the case of $G/K=SU(n+1)/S(U(1)\times U(n)),\,Sp(n+1)/(Sp(1)\times Sp(n))$ or $F_4/{\rm Spin}(9)$, we write $\triangle_+$ 
as $\triangle_+=\{\lambda,2\lambda\}$.  
The multiplicity $m_{2\lambda}$ of $2\lambda$ is given by 
\begin{equation}\label{eq:multi}
m_{2\lambda}=\left\{\begin{array}{ll}
1 & ({\rm when}\,\,G/K=SU(n+1)/S(U(1)\times U(n)))\\
3 & ({\rm when}\,\,G/K=Sp(n+1)/(Sp(1)\times Sp(n)))\\
7 & ({\rm when}\,\,G/K=F_4/{\rm Spin}(9))
\end{array}\right.
\end{equation}

Verh$\acute{\rm o}$czki \cite{V} described explicitly the eigenvalues (i.e., the principal curvatures) of the shape operators of the orbits of 
the isotropy group actions $K\curvearrowright G/K$ by using the positive restricted roots (i.e., the elements of $\triangle_+$ (see Theorem 1 of \cite{V}).  
By using the explicit descriptions of principal curvatures, (\ref{eq:root-system}) and (\ref{eq:multi}), 
we can explicitly described the above constant mean curvature $h(s)$ in the case of 
$G/K=SO(n+1)/SO(n)$, $SU(n+1)/S(U(1)\times U(n))$, $Sp(n+1)/(Sp(1)\times Sp(n))$ or $F_4/{\rm Spin}(9)$ as follows.  

\begin{lemma}
The mean curvatyre $h(s)$ of the principle orbit $r^{-1}(s)$ is given by 
\begin{equation}\label{eq:cmc}
h(s)=\left\{
\begin{array}{ll}
\displaystyle\sqrt{\frac{n-1}{2}}\cdot\frac{1}{\tan{\frac{s}{\sqrt{2(n-1)}}}} & (G/K=SO(n+1)/SO(n))\\
\displaystyle\frac{2n-1-\tan^2{\frac{s}{2\sqrt{n+1}}}}{2\sqrt{n+1}\cdot\tan{\frac{s}{2\sqrt{n+1}}}} & (G/K=SU(n+1)/S(U(1)\times U(n)))\\
\displaystyle\frac{4n-1-3\tan^2{\frac{s}{2\sqrt{2(n+2)}}}}{2\sqrt{2(n+2)}\cdot\tan{\frac{s}{2\sqrt{2(n+2)}}}} & (G/K=Sp(n+1)/(Sp(1)\times Sp(n)))\\
\displaystyle{\frac{(16-7\tan^2{as})a}{\tan{as}}} & (G/K=F_4/{\rm Spin}(9))\\
\end{array}\right.
\end{equation}
\end{lemma}

From Proposition 2.2 and Lemma 2.3, we prove Theorem 1.1.  

\vspace{0.5truecm}

\noindent
{\it Proof of Theorem 1.1.}\ \ 
We consider the case of $G/K=SU(n+1)/S(U(1)\times U(n))$.  
Define a function $\psi$ by 
$$\psi(x):=V'(2\sqrt{n+1}\arctan x)\quad\,\,(x\in r(M)).$$
From (\ref{eq:ode-cohomo1}) and (\ref{eq:cmc}), it is shown that $\psi$ satisfies the following ODE:
\begin{equation}\label{eq:ode-psi-case2}
\psi'(x)=\frac{2\sqrt{n+1}}{1+x^2}\left(1+\psi(x)^2\right)\left(1-\frac{2n-1-x^2}{2\sqrt{n+1}\,x}\psi(x)\right)\quad(x>0)
\end{equation}
We shall analyze the shape of the solution $\psi$ of (\ref{eq:ode-psi-case2}) to recognize the shape of $V$.
Define a function $\eta:[0,\infty)\setminus\{\sqrt{2n-1}\}\to\mathbb R$ by $\eta(x):=\frac{2\sqrt{n+1}x}{2n-1-x^2}$.  
From (\ref{eq:ode-psi-case2}), we can show the following facts directly:
\begin{itemize}

\item[(i)] When $x\in(0,\sqrt{2n-1})$, 
\begin{itemize}
\item[(i-a)] if $\psi(x)>\eta(x)$, then $\psi'(x)<0$
\item[(i-b)] if $\psi(x)<\eta(x)$, then $\psi'(x)>0$
\item[(i-c)] if $\psi(x)=\eta(x)$, then $\psi'(x)=0;$
\end{itemize}
\item[(ii)] When $x\in(\sqrt{2n-1},\infty)$, 
\begin{itemize}
\item[(ii-a)] if $\psi(x)>\eta(x)$, then $\psi'(x)>0$
\item[(ii-b)] if $\psi(x)<\eta(x)$, then $\psi'(x)<0$
\item[(ii-c)] if $\psi(x)=\eta(x)$, then $\psi'(x)=0$
\end{itemize}
\item[(iii)] When $x=\sqrt{2n-1}$, $\psi'(x)>0$. 
\end{itemize}

\begin{figure}[H]
\centering
\scalebox{1.0}{\input{eta-graph.tex}}
\caption{The graph of $\eta$}
\end{figure}

\noindent
Next we shall show that the following fact for $\psi$ holds.  

\vspace{0.15truecm}

($\ast_1$)\ \ {\sl 
If there exists $x_0\in(0,\sqrt{2n-1})$ with $\psi(x_0)>\eta(x_0)$, $\lim\limits_{x\downarrow x_1}{\psi(x)}=\infty$ holds for some 

$x_1\in(0,x_0)$ (see Figure 2.2).
}  

\vspace{0.15truecm}

\noindent
Take any $x\in(0,x_0)$.  Then, by using (i-a), we can show 
\begin{align*}
\psi'(x)&=\frac{2\sqrt{n+1}}{1+x^2}\left(1+\psi(x)^2\right)\left(1-\frac{\psi(x)}{\eta(x)}\right)\\
&<\frac{2\sqrt{n+1}}{1+x^2}\left(1+\psi(x)^2\right)\left(1-\frac{\psi(x_0)}{\eta(x)}\right).
\end{align*}
Therefore, we have 
\begin{equation*}
\frac{\psi'(x)}{1+\psi(x)^2}<\frac{2\sqrt{n+1}}{1+x^2}-\frac{(2n-1)\psi(x_0)}{x(1+x^2)}+\frac{\psi(x_0)x}{1+x^2}.
\end{equation*}
By integrating both-hand sides of this inequality from $x$ to $x_0$, we obtain
\begin{align*}
\arctan{\psi(x_0)}-\arctan{\psi(x)}<&-2\sqrt{n+1}\arctan{x}+(2n-1)\psi(x_0)\log{\frac{x}{\sqrt{1+x^2}}}\\
&-\frac{\psi(x_0)}{2}\log{(1+x^2)}+2\sqrt{n+1}\arctan{x_0}\\
&-(2n-1)\psi(x_0)\log{\frac{x_0}{\sqrt{1+{x_0}^2}}}+\frac{\psi(x_0)}{2}\log{(1+{x_0}^2)}=:h_1(x).
\end{align*}
and hence 
$$\psi(x)>\tan{(-h_1(x)+\arctan{\psi(x_0)})}.$$
On the other hand, $h_1$ is increasing on $(0,x_0)$ and the following relations hold:
$$h_1(x_0)=0\quad\,\,{\rm and}\quad\,\,\lim_{x\downarrow0}{h_1(x)}=-\infty.$$
Therefore, there exists $\bar{x}_1\in(0,x_0)$ such that
$$\psi(x)>\tan{(-h_1(x)+\arctan{\psi(x_0)})}\,\,\rightarrow\,\,\infty\quad(x\rightarrow \bar{x}_1).$$
Thus the fact ($\ast_1$) is shown.  

\begin{figure}[H]
\centering
\scalebox{1.0}{\input{psi2-lemma1.tex}}
\caption{The behavior I of the graph of $\psi$}
\end{figure}

\noindent
Similarly, by using (i-b), we can show the following facts for $\psi$;

\vspace{0.15truecm}

($\ast_2$)\ \ {\sl If there exists $x_0\in(0,\sqrt{2n-1})$ with $\psi(x_0)<0$, $\lim\limits_{x\downarrow x_1}{\psi(x)}=-\infty$ holds for some 

$x_1\in(0,x_0)$ (see Figure 2.3);}  

\vspace{0.15truecm}

\noindent
Also, by using (ii-a), we can show the following facts for $\psi$;

\vspace{0.15truecm}

($\ast_3$)\ \ If there exists $x_0\in(\sqrt{2n-1},\infty)$ with $\psi(x_0)>0$, $\lim\limits_{x\uparrow x_1}{\psi(x)}=\infty$ holds for some 

$x_1\in(x_0,\infty)$ (see Figure 2.4);

\vspace{0.15truecm}

\noindent
Also, by using (ii-b), we can show the following facts for $\psi$;

\vspace{0.15truecm}

($\ast_4$)\ \ {\sl 
If there exists $x_0\in(\sqrt{2n-1},\infty)$ with $\psi(x_0)<\eta(x_0)$, $\lim\limits_{x\uparrow x_1}{\psi(x)}=-\infty$ holds for some 

$x_1\in(x_0,\infty)$ (see Figure 2.5).  
}

\begin{figure}[H]
\centering
\scalebox{1.0}{\input{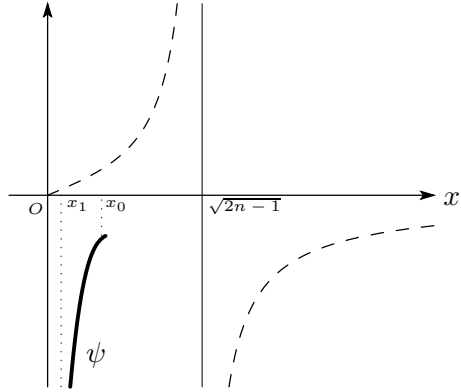}}
\caption{The behavior II of the graph of $\psi$}
\end{figure}

\begin{figure}[H]
\centering
\scalebox{1.0}{\input{psi2-lemma3.tex}}
\caption{The behavior III of the graph of $\psi$}
\end{figure}

\begin{figure}[H]
\centering
\scalebox{1.0}{\input{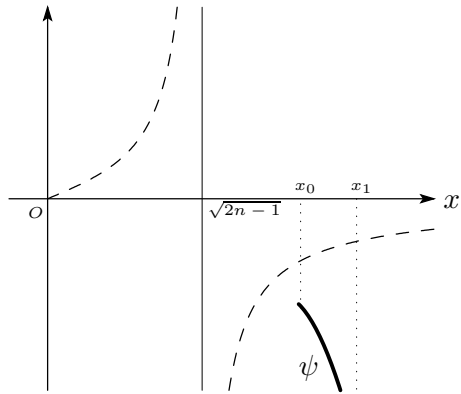}}
\caption{The behavior IV of the graph of $\psi$}
\end{figure}

\noindent
From 
the facts $(\ast_1)-(\ast_4)$, the shape of the graph of the solution $\psi$ is one of the curves as in Figures \ref{psi-TypeI}-\ref{psi-TypeV} 
in the case of $G/K=SU/S(U(1)\times U(n))$.  
\begin{figure}[H]
\centering
\scalebox{0.9}{\input{psi2-shape1.tex}}
\caption{Type I}
\label{psi-TypeI}
\end{figure}
\begin{figure}[H]
\centering
\begin{minipage}{0.45\columnwidth}
\centering
\scalebox{0.9}{\input{psi2-shape2.tex}}
\caption{Type II}
\end{minipage}
\begin{minipage}{0.45\columnwidth}
\centering
\scalebox{0.9}{\input{psi2-shape3.tex}}
\caption{Type III}
\end{minipage}
\end{figure}
\begin{figure}[h]
\centering
\begin{minipage}{0.45\columnwidth}
\centering
\scalebox{0.9}{\input{psi2-shape4.tex}}
\caption{Type IV}
\end{minipage}
\begin{minipage}{0.45\columnwidth}
\centering
\scalebox{0.9}{\input{psi2-shape5.tex}}
\caption{Type V}
\label{psi-TypeV}
\end{minipage}
\end{figure}

\noindent
From (\ref{eq:cmc}), we see that the domain $V$ is included by $(0,\sqrt{n+1}\pi)$, that is, the value $\alpha$ in the statement of Theorem 1.1 is equal to 
$\sqrt{n+1}\pi$.
Hence, from the above classification of the shape of $\psi$, we can classify the shape of $V$ as in Theorem 1.1.

Similarly, in the case where $G/K=SO(n+1)/SO(n)$, $Sp(n+1)/(Sp(1)\times Sp(n))$ and $F_4/{\rm Spin}(9)$, we can classify the shape of the graph of $V$.  
In these cases, by (\ref{eq:cmc}), the value $\alpha$ in the statement of Theorem 1.1 is given by
$$\alpha=\left\{
\begin{array}{ll}
\displaystyle{\sqrt{\frac{n-1}{2}}\pi} & ({\rm when}\,\,G/K=SO(n+1)/SO(n))\\
\displaystyle{\sqrt{2(n+2)}\pi} & ({\rm when}\,\,G/K=Sp(n+1)/(Sp(1)\times Sp(n)))\\
\displaystyle{\frac{\pi}{4a}} & ({\rm when}\,\,G/K=F_4/{\rm Spin}(9)).
\end{array}
\right.$$
\qed

\section{The case of cohomogenity two}
In this section, we consider the case where $G/K$ is a higher rank irreducible symmetric space of compact type and a translator given  by a graph of a function $u$ on 
a $H$-invariant domian $M$ of $G/K$ which is invariant under a Hermann action $H\curvearrowright G/K$ of cohomegeneity two.  
Assume that $H\curvearrowright G/K$ is commutative, that is, $\theta_K\circ\theta_H=\theta_H\circ\theta_K$ holds for the involutions $\theta_K$ and $\theta_H$ of $G$ satisfying 
$(\textrm{Fix}\theta_K)_0\subset K\subset\textrm{Fix}\theta_K$ and $(\textrm{Fix}\theta_H)_0\subset H\subset\textrm{Fix}\theta_H$.
Let $r=(r_1,r_2):G/K\to\mathbb{R}^2$ be a map on $G/K$ with $g(\nabla r_i,\nabla r_j)=\delta_{ij}$ ($i,j\in\{1,2\}$) whose level sets give the orbits of 
the action $H\curvearrowright G/K$.  Then, a function $u$ is described as $u=V\circ r$ for some function $V$ on $r(M)$.  

By using Theorem 1.1, we prove Theorem 1.2.  

\vspace{0.5truecm}

\noindent
{\it Proof of Theorem 3.1.2.}\ \ 
The function $V$ is regarded as a function on $\mathfrak a$ through the Euclidean coordinate $(x_1,x_2):\mathfrak a\to\mathbb R^2$ of $\mathfrak a$.  
As stated in  Introduction, we may assume that 
\begin{equation}\label{eq:r}
r|_{\exp_o(\mathcal C)}=-(x_1,x_2)\circ(\exp_o|_{\mathcal C})^{-1}.
\end{equation}
From $u=V\circ r$, we find
\begin{align*}
&\nabla u=\sum_{i=1}^2\frac{\partial V}{\partial x_i}\nabla r_i,\quad\|\nabla u\|^2=|\nabla V|^2,\\
&{\rm div}\left(\frac{\nabla u}{\sqrt{1+\|\nabla u\|^2}}\right)=\sum_{1=1}^2\Delta r_i\frac{\partial V}{\partial x_i}+\Delta V-\frac{1}{1+|\nabla V|^2}\sum_{i,j=1}^2\frac{\partial V}{\partial x_i}\frac{\partial V}{\partial x_j}\frac{\partial^2 V}{\partial x_i\partial x_j}.
\end{align*}
From these relations, we can show that the PDE (\ref{eq:graph-soliton}) is reduced to 
\begin{equation}\label{eq:pde-cohomo2-2}
\sum_{i,j=1}^2\frac{\partial V}{\partial x_i}\frac{\partial V}{\partial x_j}\frac{\partial^2 V}{\partial x_i\partial x_j}-\left(1+|\nabla V|^2\right)\left(\sum_{i=1}^2
\Delta r_i\frac{\partial V}{\partial x_i}+\Delta V-1\right)=0.
\end{equation}
Since the tangent vector $\vX$ on $\exp_o(\mathcal C)$ is defined by assigning the mean curvature vector of the orbit $H\cdot w$ at $w$ to 
each $w\in\exp_o(\mathcal C)$, we have
\begin{equation*}
\langle\vX,\nabla r_i\rangle=-\Delta r_i.
\end{equation*}
From this relation and (\ref{eq:r}), we find
\begin{align*}
\bm X=\sum_{i=1}^2\Delta r_i\left(-\nabla r_i\right).
\end{align*}
Therefore, from (\ref{eq:r}), we obtain $X_i=\Delta r_i$.
\qed

By using Theorem 1.2, we prove Theorem 1.3.  

\vspace{0.5truecm}

\noindent{\it Proof of Theorem 1.3.}\ \ 
Assume that $V$ is constant along each level set of $\rho$, where $\rho$ is the convex function with $\nabla\rho=\vX$.
Then, $\nabla V$ is described as $\nabla V=F\vX$ for some function $F$ on $\mathcal{C}$.  
Clearly we have 
\begin{equation*}
\frac{\partial V}{\partial x_i}=FX_i,\quad\frac{\partial^2 V}{\partial x_i\partial x_j}=\frac{\partial F}{\partial x_i}X_j+F\frac{\partial X_j}{\partial x_i}\quad(i,j\in\{1,2\}).
\end{equation*}
Also, we can derive 
\begin{align*}
|&\nabla V|^2=|\vX|^2F^2,\quad\Delta V=\vX(F)+(\textrm{div}\vX)F,\quad\sum_{i=1}^2X_i\frac{\partial V}{\partial x_i}=|\vX|^2F,\\
&\sum_{i,j=1}^2\frac{\partial^2V}{\partial x_i\partial x_j}\frac{\partial V}{\partial x_i}\frac{\partial V}{\partial x_j}=|\vX|^2F^2\vX(F)+\frac{1}{2}\vX(|\vX|^2)F^3.
\end{align*}
Hence we obtain
\begin{align*}
&\sum_{i,j=1}^2\frac{\partial^2V}{\partial x_i\partial x_j}\frac{\partial V}{\partial x_i}\frac{\partial V}{\partial x_j}
-\left(1+|\nabla V|^2\right)\left(\sum_{i=1}^2X_i\frac{\partial V}{\partial x_i}+\Delta V-1\right)\\
=&|\vX|^2F^2\vX(F)+\frac{1}{2}\vX(|\vX|^2)F^3-\left(1+|\vX|^2F^2\right)\left(|\vX|^2F+\vX(F)+(\textrm{div}\vX)F-1\right)\\
=&-\vX(F)+\frac{1}{2}\vX(|\vX|^2)F^3-\left(1+|\vX|^2F^2\right)\left(\left(|\vX|^2+(\textrm{div}\vX)\right)F-1\right).
\end{align*}
Therefore, from Theorem 1.2, we can derive Theorem 1.3.  \qed

\vspace{0.5truecm}

Denote by $(a_1,a_2)$ be the minimum point of $\rho$.  Let $F$ be a solution of the partial differential equation (\ref{eq:pde2-cohomo2}) and 
$c:(-\infty,t_0)\to\mathfrak a$ be an integral curve of $\vX$, where we note that $\lim\limits_{t\to-\infty}c(t)=(a_1,a_2)$ holds.  
Set $\widehat{F}_c:=F\circ c$.
Then $\widehat{F}_c$ satisfies
\begin{equation}\label{eq:pde-Fhat}
\widehat{F}_c'(t)=\langle c''(t),c'(t)\rangle\widehat{F}_{c}(t)^3-(1+|c'(t)|^2\widehat{F}_{c}(t)^2)((|c'(t)|^2+(\textrm{div}\vX)_{c(t)})\widehat{F}_{c}(t)-1).
\end{equation}
From $\nabla V=F\vX$, we have
\begin{equation*}
(V\circ c)'(t)=\widehat F_c(t)|\vX_{c(t)}|^2
\end{equation*} 
and hence 
\begin{equation}\label{eq:Vc}
(V\circ c)(t):=\int_{t_{\ast}}^t\widehat F_c(\tau)|\vX_{c(\tau)}|^2d\tau+(V\circ c)(t_{\ast})),
\end{equation}
where $t_{\ast}$ is any element of $(-\infty,t_0)$.  Thus we can calculate the function $V$ from the data of $F$. 

We shall consider the case where the Hermann action $H\curvearrowright G/K$ is the dual action of the Hermann type action $SO_0(1,2)\curvearrowright SL(3,\mathbb R)/SO(3)$.  
This action corresponds to $\rho_1(SO(3))\curvearrowright SU(3)/SO(3)$ in Table 1 of \cite{K}.  
In this case, there exists an integral curve $c:(-\infty,t_0)\to\mathfrak a$ of $\vX$ satisfying $x_2\circ c=0$.
The component $X_i$ of $\vX$ and the domain $\mathcal C$ are given as in the following table.  
\begin{table}[H]
\centering
\begin{tabular}{|c|c|}
\hline
$H\curvearrowright G/K$&{\scriptsize$\rho_1(SO(3))\curvearrowright SU(3)/SO(3)$}\\
\hline
$X_1$&{\scriptsize$\textrm{tan}(x_1+\sqrt{3}x_2)-2\textrm{cot}(2x_1)+\textrm{tan}(x_1-\sqrt{3}x_2)$}\\
\hline
$X_2$&{\scriptsize$\sqrt{3}\textrm{tan}(x_1+\sqrt{3}x_2)-\sqrt{3}\textrm{tan}(x_1-\sqrt{3}x_2)$}\\
\hline
$\mathcal C$&{\scriptsize$\{(x_1,x_2)~|~x_1>0,~x_2>\frac{1}{\sqrt{3}}x_1-\frac{\pi}{2\sqrt{3}},~x_2<-\frac{1}{\sqrt{3}}x_1+\frac{\pi}{2\sqrt{3}}\}$}\\
\hline
\end{tabular}
\caption{The datas of $X_1,X_2$ and $\mathcal C$ in $\rho_1(SO(3))\curvearrowright SU(3)/SO(3)$-case}
\end{table}

\noindent 
In this case, $(a_1,a_2)$ is equal to $(\frac{\pi}{6},0)$ and there exists an integral curve $c:(-\infty,t_0)\to\mathfrak a$ satisfying 
$\lim\limits_{t\to-\infty}c(t)=(\frac{\pi}{6},0),\,\,\,\lim\limits_{t\to t_0}c(t)=(\frac{\pi}{2},0)$ and $x_2\circ c=0$.  
For the simplicity, set $c_1:=x_1\circ c$.  Also, set $\widetilde{F}:=|c'|\widehat{F}_c$.
Take $t_{\ast}\in(-\infty,t_0)$.  Since $\vX_{c(t)}=c'(t)=(c_1'(t),0)$ and hence $\|\vX_{c(t)}\|=c_1'(t)$, we have 
\begin{align*}
V(x_1,0)&=V(c_1(t),0)=(V\circ c)(t)\\
&=\int_{t_{\ast}}^t\widehat F_c(\tau)\|\vX_{c(\tau)}\|^2d\tau+V(c_1(t_{\ast}),0)\\
&=\int_{t_{\ast}}^t\widehat F_c(\tau)c_1'(\tau)^2d\tau+V(c_1(t_{\ast}),0)\\
&=\int_{t_{\ast}}^t\widetilde F(\tau)c_1'(\tau)\,d\tau+V(c_1(t_{\ast}),0)\\
&=\int_{c_1(t_{\ast})}^{x_1}(\widetilde F\circ c_1^{-1})(x)\,dx+V(c_1(t_{\ast}),0)\quad\left(\frac{\pi}{6}\leq x_1\leq\frac{\pi}{2}\right)
\end{align*}
The shape of the graph of the restriction $V|_{[\frac{\pi}{6},\frac{\pi}{2}]\times\{0\}}$ of $V$ to 
$\overline{c((-\infty,t_0))}=\left[\frac{\pi}{6},\frac{\pi}{2}\right]\times\{0\}$ is dominated by $\widetilde F\circ c_1^{-1}$.  
So, we shall investigate the shape of the graph of $\widetilde F\circ c_1^{-1}$.  
According to (\ref{eq:pde-Fhat}), $\widetilde{F}$ satisfies
\begin{equation}\label{eq:tildeF}
\widetilde{F}'(t)=\left(\widetilde{F}(t)^2+1\right)\left(\left(\frac{\langle c'(t),c''(t)\rangle}{|c'(t)|^2}-\left(|c'(t)|^2+(\textrm{div}\vX)_{c(t)}\right)\right)\widetilde{F}(t)+|c'(t)|\right),
\end{equation}
On the othe hand, we have 
\begin{align}\label{1st-relation}
|c'(t)|^2=9\textrm{tan}^2(c_1(t))+\textrm{cot}^2(c_1(t))-6,
\end{align}
\begin{align}\label{2nd-relation}
\langle c'(t),c''(t)\rangle&=\frac{27\textrm{tan}^2(c_1(t))}{\textrm{cos}^2(c_1(t))}+\frac{\textrm{cot}^2(c_1(t))}{\textrm{sin}^2(c_1(t))}-\frac{9}{\textrm{cos}^2(c_1(t))}-\frac{3}{\textrm{sin}^2(c_1(t))}
\end{align}
and 
\begin{align}\label{3rd-relation}
(\textrm{div}\vX)_{c(t)}=\frac{9}{\textrm{cos}^2(c_1(t))}+\frac{1}{\textrm{sin}^2(c_1(t))}.
\end{align}
Define $\widetilde{\eta}(t)$ by 
\begin{equation*}
\widetilde{\eta}(t):=-\frac{|c'(t)|^3}{\langle c'(t),c''(t)\rangle-|c'(t)|^2\left(|c'(t)|^2+(\textrm{div}\vX)_{c(t)}\right)}.
\end{equation*}
From (\ref{1st-relation}), (\ref{2nd-relation}) and (\ref{3rd-relation}), we have 
\begin{equation*}
\widetilde{\eta}(t)=\frac{3\textrm{tan}(c_1(t))-\frac{1}{\textrm{tan}(c_1(t))}}{15\textrm{tan}^2(c_1(t))+\frac{1}{\textrm{tan}^2(c_1(t))}}
\end{equation*}
Set $z(t):=\textrm{tan}(c_1(t))$ and define $\sigma(z)$ by
\begin{equation*}
\sigma(z):=\frac{3z-\frac{1}{z}}{15z^2+\frac{1}{z^2}}.
\end{equation*}
From $\widetilde{\eta}=\sigma(z)$, we find 
\begin{equation}\label{eq:deta}
\widetilde{\eta}'(t)=\left(1+z(t)^2\right)\left(3z(t)-\frac{1}{z(t)}\right)\sigma'(z(t)).
\end{equation}
Since $\sigma$ satisfies
\begin{equation*}
\sigma'(z)=\frac{-45z^6+45z^4+9z^2-1}{z^4\left(15z^2+\frac{1}{z^2}\right)^2},
\end{equation*}
there exists $z_0\in(\frac{1}{\sqrt{3}},\infty)$ with $\sigma'(z_0)=0$ satisfying $\sigma'(z)>0$ for all $z\in(\frac{1}{\sqrt{3}},z_0)$ and $\sigma'(z)<0$ for all $z\in(z_0,\infty)$.
Define $t_1\in(-\infty,t_0)$ by $\textrm{tan}(c_1(t_1))=z_0$.
Then, from (\ref{eq:deta}), we find $\widetilde{\eta}'(t_1)=0$, $\widetilde{\eta}'(t)>0$ for all $t\in(-\infty,t_1)$ and $\widetilde{\eta}'(t)<0$ for all $t\in(t_1,t_0)$.
Therefore, the graph of $\widetilde{\eta}\circ c_1^{-1}$ is as in Figure 3.1. 
\begin{figure}[H]
\centering
\scalebox{1.0}{
{\unitlength 0.1in%
\begin{picture}(20.3000,12.0000)(2.0000,-14.0000)%
\special{pn 13}%
\special{pn 13}%
\special{pa 400 1196}%
\special{pa 405 1190}%
\special{pa 410 1181}%
\special{pa 420 1159}%
\special{pa 425 1146}%
\special{pa 427 1140}%
\special{fp}%
\special{pa 448 1080}%
\special{pa 455 1056}%
\special{pa 460 1039}%
\special{pa 465 1023}%
\special{pa 466 1020}%
\special{fp}%
\special{pa 484 959}%
\special{pa 485 956}%
\special{pa 503 900}%
\special{fp}%
\special{pa 523 839}%
\special{pa 525 832}%
\special{pa 535 804}%
\special{pa 544 781}%
\special{fp}%
\special{pa 569 722}%
\special{pa 570 720}%
\special{pa 575 709}%
\special{pa 585 689}%
\special{pa 595 671}%
\special{pa 597 667}%
\special{fp}%
\special{pa 634 615}%
\special{pa 640 607}%
\special{pa 660 587}%
\special{pa 679 572}%
\special{fp}%
\special{pa 736 546}%
\special{pa 755 542}%
\special{pa 760 542}%
\special{pa 765 541}%
\special{pa 795 541}%
\special{pa 798 542}%
\special{fp}%
\special{pa 861 552}%
\special{pa 865 554}%
\special{pa 870 555}%
\special{pa 875 557}%
\special{pa 880 558}%
\special{pa 890 562}%
\special{pa 895 563}%
\special{pa 919 573}%
\special{fp}%
\special{pa 978 598}%
\special{pa 985 601}%
\special{pa 990 604}%
\special{pa 995 606}%
\special{pa 1000 609}%
\special{pa 1005 611}%
\special{pa 1010 614}%
\special{pa 1015 616}%
\special{pa 1020 619}%
\special{pa 1025 621}%
\special{pa 1033 626}%
\special{fp}%
\special{pa 1089 656}%
\special{pa 1105 665}%
\special{pa 1110 667}%
\special{pa 1135 682}%
\special{pa 1140 684}%
\special{pa 1143 686}%
\special{fp}%
\special{pa 1198 718}%
\special{pa 1252 750}%
\special{fp}%
\special{pa 1307 782}%
\special{pa 1360 814}%
\special{fp}%
\special{pa 1415 847}%
\special{pa 1468 879}%
\special{fp}%
\special{pa 1523 912}%
\special{pa 1576 944}%
\special{fp}%
\special{pa 1631 977}%
\special{pa 1645 985}%
\special{pa 1650 987}%
\special{pa 1685 1008}%
\special{fp}%
\special{pa 1739 1041}%
\special{pa 1770 1059}%
\special{pa 1775 1061}%
\special{pa 1793 1072}%
\special{fp}%
\special{pa 1848 1105}%
\special{pa 1865 1115}%
\special{pa 1870 1117}%
\special{pa 1902 1136}%
\special{fp}%
\special{pa 1956 1169}%
\special{pa 1960 1170}%
\special{pa 2010 1200}%
\special{fp}%
\special{pn 13}%
\special{pa 2024 1192}%
\special{pa 2040 1182}%
\special{pa 2045 1180}%
\special{pa 2058 1172}%
\special{ip}%
\special{pa 2071 1164}%
\special{pa 2105 1144}%
\special{ip}%
\special{pa 2119 1136}%
\special{pa 2125 1132}%
\special{pa 2130 1130}%
\special{pa 2153 1116}%
\special{ip}%
\special{pa 2166 1108}%
\special{pa 2200 1088}%
\special{ip}%
%
\special{pn 8}%
\special{pa 200 1200}%
\special{pa 2200 1200}%
\special{fp}%
\special{sh 1}%
\special{pa 2200 1200}%
\special{pa 2133 1180}%
\special{pa 2147 1200}%
\special{pa 2133 1220}%
\special{pa 2200 1200}%
\special{fp}%
%
\special{pn 4}%
\special{pa 400 200}%
\special{pa 400 1400}%
\special{dt 0.045}%
%
\special{pn 4}%
\special{pa 2000 200}%
\special{pa 2000 1400}%
\special{dt 0.045}%
\put(22.3000,-12.4000){\makebox(0,0)[lb]{{\tiny$x_1$}}}%
\put(2.8000,-13.3000){\makebox(0,0)[lb]{{\tiny$\frac{\pi}{6}$}}}%
\put(20.3000,-13.3000){\makebox(0,0)[lb]{{\tiny$\frac{\pi}{2}$}}}%
\put(10.5000,-5.5000){\makebox(0,0)[lb]{$\widetilde{\eta}\circ {c_1}^{-1}$}}%
%
\special{pn 4}%
\special{pa 760 540}%
\special{pa 760 1200}%
\special{dt 0.045}%
\put(6.3000,-13.3000){\makebox(0,0)[lb]{{\tiny$c_1(t_1)$}}}%
\end{picture}}%
}
\caption{The graph of $\widetilde{\eta}\circ c_1^{-1}$}
\end{figure}
Here, from $\eta(t)>0$, we find $\frac{\langle c'(t),c''(t)\rangle}{|c'(t)|^2}-\left(|c'(t)|^2+(\textrm{div}\vX)_{c(t)}\right)<0$.
According to (\ref{eq:tildeF}), $\widetilde{F}'(t)=0$ if and only if $\widetilde{F}(t)=\widetilde{\eta}(t)$.
Also, $\widetilde{F}'(t)>0$ if and only if $\widetilde{F}(t)<\widetilde{\eta}(t)$, and $\widetilde{F}'(t)<0$ if and only if $\widetilde{F}(t)>\widetilde{\eta}(t)$.

Next we show the behavior of the function $\widetilde{F}\circ c_1^{-1}$ on the both sides of the domain of $\widetilde{F}\circ c_1^{-1}$.
From (\ref{eq:tildeF}), $\widetilde{F}\circ c_1^{-1}$ satisfies
\begin{equation}\label{eq:tildeFc}
(\widetilde{F}\circ c_1^{-1})'(x)=\left((\widetilde{F}\circ c_1^{-1})(x)^2+1\right)\left(1-\left(5\tan{x}-\frac{1}{\tan{x}}+\frac{8\tan{x}}{3\textrm{tan}^2x-1}\right)(\widetilde{F}\circ c_1^{-1})(x)\right).
\end{equation}
Then we shall show that the following fact for the behavior of $\widetilde{F}\circ c_1^{-1}$ near $x_1=\frac{\pi}{2}$ holds.

\vspace{0.15truecm}

($\ast_5$)\ \ {\sl 
If there exists $x_0\in(c_1(t_1),\frac{\pi}{2})$ such that $(\widetilde{F}\circ c_1^{-1})(x_0)>\widetilde{\eta}(x_0)$, $\lim\limits_{x\uparrow \frac{\pi}{2}}{(\widetilde{F}\circ c_1^{-1})(x)}=0$ holds.
}  

\vspace{0.15truecm}

\noindent
Assume that there exists a positive constant $M>0$ with $(\widetilde{F}\circ c_1^{-1})(x)>M$ for all $x\in(x_0,\frac{\pi}{2})$.
Take any $x\in(x_0,\frac{\pi}{2})$.  Then, from $(\widetilde{F}\circ c_1^{-1})(x)<0$, we can show 
\begin{align*}
\frac{(\widetilde{F}\circ c_1^{-1})'(x)}{(\widetilde{F}\circ c_1^{-1})(x)^2+1}&=1-\left(5\tan{x}-\frac{1}{\tan{x}}+\frac{8\tan{x}}{3\textrm{tan}^2x-1}\right)(\widetilde{F}\circ c_1^{-1})(x)\\
&<1-\left(5\tan{x}-\frac{1}{\tan{x}}+\frac{8\tan{x}}{3\textrm{tan}^2x-1}\right)M.
\end{align*}
By integrating both-hand sides of this inequality from $x_0$ to $x$, we obtain
\begin{align*}
\arctan{(\widetilde{F}\circ c_1^{-1})(x)}-&\arctan{(\widetilde{F}\circ c_1^{-1})(x_0)}\\
<&~x+M(5\log(\cos{x})+\log{(\sin{x})}-\log{(1-2\cos{(2x)})})\\
&-x_0-M(5\log(\cos{x_0})+\log{(\sin{x_0})}-\log{(1-2\cos{(2x_0)})})=:h_2(x).
\end{align*}
and hence 
$$(\widetilde{F}\circ c_1^{-1})(x)<\tan{(h_2(x)+\arctan{(\widetilde{F}\circ c_1^{-1})(x_0)})}.$$
On the other hand, $h_2$ is decreasing on $(x_0,\frac{\pi}{2})$ and the following relations hold:
$$h_2(x_0)=0\quad\,\,{\rm and}\quad\,\,\lim_{x\uparrow\frac{\pi}{2}}{h_2(x)}=-\infty.$$
Therefore, there exists $\bar{x}_1\in(x_0,\frac{\pi}{2})$ such that
$$(\widetilde{F}\circ c_1^{-1})(x)<\tan{(h_2(x)+\arctan{(\widetilde{F}\circ c_1^{-1})(x_0)})}\,\,\rightarrow\,\,-\infty\quad(x\rightarrow \bar{x}_1).$$
Then $(\widetilde{F}\circ c_1^{-1})(x)=0<M$ for some $x\in(x_0,\bar{x}_1)$.
This is a contradiction.
Thus the fact ($\ast_5$) is shown.

Also, we shall show the following fact for the behavior of $\widetilde{F}\circ c_1^{-1}$ near some point $x_1=x_0\in(\frac{\pi}{6},c_1(t_1))$.

\vspace{0.15truecm}

($\ast_6$)\ \ {\sl 
If there exists $x_0\in(\frac{\pi}{6},c_1(t_1))$ such that $(\widetilde{F}\circ c_1^{-1})(x_0)>\widetilde{\eta}(x_0)$, $\lim\limits_{x\downarrow \bar{x}_0}{(\widetilde{F}\circ c_1^{-1})(x)}=\infty$ holds for some $\bar{x}_0\in(\frac{\pi}{6},x_0)$.
}  

\vspace{0.15truecm}

\noindent
Take any $x\in(\frac{\pi}{6},x_0)$.
Then, from $(\widetilde{F}\circ c_1^{-1})(x)<0$, we can show 
\begin{align*}
\frac{(\widetilde{F}\circ c_1^{-1})'(x)}{(\widetilde{F}\circ c_1^{-1})(x)^2+1}&=1-\left(5\tan{x}-\frac{1}{\tan{x}}+\frac{8\tan{x}}{3\textrm{tan}^2x-1}\right)(\widetilde{F}\circ c_1^{-1})(x)\\
&<1-\left(5\tan{x}-\frac{1}{\tan{x}}+\frac{8\tan{x}}{3\textrm{tan}^2x-1}\right)(\widetilde{F}\circ c_1^{-1})(x_0).
\end{align*}
By integrating both-hand sides of this inequality from $x$ to $x_0$, we obtain
\begin{align*}
\arctan{(\widetilde{F}\circ c_1^{-1})(x)}-&\arctan{(\widetilde{F}\circ c_1^{-1})(x_0)}\\
>&~x+((\widetilde{F}\circ c_1^{-1})(x_0))(5\log(\cos{x})+\log{(\sin{x})}-\log{(1-2\cos{(2x)})})\\
&-x_0-((\widetilde{F}\circ c_1^{-1})(x_0))(5\log(\cos{x_0})+\log{(\sin{x_0})}-\log{(1-2\cos{(2x_0)})})=:h_3(x).
\end{align*}
and hence 
$$(\widetilde{F}\circ c_1^{-1})(x)>\tan{(h_3(x)+\arctan{(\widetilde{F}\circ c_1^{-1})(x_0)})}.$$
On the other hand, $h_3$ is decreasing on $(\frac{\pi}{6},x_0)$ and the following relations hold:
$$h_3(x_0)=0\quad\,\,{\rm and}\quad\,\,\lim_{x\downarrow\frac{\pi}{6}}{h_3(x)}=\infty.$$
Therefore, there exists $\bar{x}_1\in(\frac{\pi}{6},x_0)$ such that
$$(\widetilde{F}\circ c_1^{-1})(x)>\tan{(h_3(x)+\arctan{(\widetilde{F}\circ c_1^{-1})(x_0)})}\,\,\rightarrow\,\,\infty\quad(x\rightarrow \bar{x}_1).$$
Thus the fact ($\ast_6$) is shown.

Similarly, we can show the following fact for the behavior of $\widetilde{F}\circ c_1^{-1}$ near  some point $x_1=x_0\in(\frac{\pi}{6},c_1(t_1))$.
Here, note the fact that $(\widetilde{F}\circ c_1^{-1})(x)>0$ when $(\widetilde{F}\circ c_1^{-1})(x)<(\widetilde{\eta}\circ c_1^{-1})(x)$ for all $x\in(\frac{\pi}{6},\frac{\pi}{2})$

\vspace{0.15truecm}

($\ast_7$)\ \ {\sl 
If there exists $x_0\in(\frac{\pi}{6},c_1(t_1))$ such that $(\widetilde{F}\circ c_1^{-1})(x_0)<0$, $\lim\limits_{x\downarrow \bar{x}_0}{(\widetilde{F}\circ c_1^{-1})(x)}=-\infty$ holds for some $\bar{x}_0\in(\frac{\pi}{6},x_0)$.
}  

\vspace{0.15truecm}

From the facts $(\ast_5)-(\ast_7)$, the graph of $\widetilde{F}\circ c_1^{-1}$ is as in one of Figures \ref{tilde-F-graph1}-\ref{tilde-F-graph5}.  
\vspace{0.5cm}
\begin{figure}[H]
\centering
\scalebox{1.0}{
{\unitlength 0.1in%
\begin{picture}(20.3000,20.0000)(2.0000,-22.0000)%
\special{pn 8}%
\special{pn 8}%
\special{pa 400 1196}%
\special{pa 405 1190}%
\special{pa 410 1181}%
\special{pa 420 1159}%
\special{pa 425 1146}%
\special{pa 427 1140}%
\special{fp}%
\special{pa 448 1080}%
\special{pa 455 1056}%
\special{pa 460 1039}%
\special{pa 465 1023}%
\special{pa 466 1020}%
\special{fp}%
\special{pa 484 959}%
\special{pa 485 956}%
\special{pa 503 900}%
\special{fp}%
\special{pa 523 839}%
\special{pa 525 832}%
\special{pa 535 804}%
\special{pa 544 781}%
\special{fp}%
\special{pa 569 722}%
\special{pa 570 720}%
\special{pa 575 709}%
\special{pa 585 689}%
\special{pa 595 671}%
\special{pa 597 667}%
\special{fp}%
\special{pa 634 615}%
\special{pa 640 607}%
\special{pa 660 587}%
\special{pa 679 572}%
\special{fp}%
\special{pa 736 546}%
\special{pa 755 542}%
\special{pa 760 542}%
\special{pa 765 541}%
\special{pa 795 541}%
\special{pa 798 542}%
\special{fp}%
\special{pa 861 552}%
\special{pa 865 554}%
\special{pa 870 555}%
\special{pa 875 557}%
\special{pa 880 558}%
\special{pa 890 562}%
\special{pa 895 563}%
\special{pa 919 573}%
\special{fp}%
\special{pa 978 598}%
\special{pa 985 601}%
\special{pa 990 604}%
\special{pa 995 606}%
\special{pa 1000 609}%
\special{pa 1005 611}%
\special{pa 1010 614}%
\special{pa 1015 616}%
\special{pa 1020 619}%
\special{pa 1025 621}%
\special{pa 1033 626}%
\special{fp}%
\special{pa 1089 656}%
\special{pa 1105 665}%
\special{pa 1110 667}%
\special{pa 1135 682}%
\special{pa 1140 684}%
\special{pa 1143 686}%
\special{fp}%
\special{pa 1198 718}%
\special{pa 1252 750}%
\special{fp}%
\special{pa 1307 782}%
\special{pa 1360 814}%
\special{fp}%
\special{pa 1415 847}%
\special{pa 1468 879}%
\special{fp}%
\special{pa 1523 912}%
\special{pa 1576 944}%
\special{fp}%
\special{pa 1631 977}%
\special{pa 1645 985}%
\special{pa 1650 987}%
\special{pa 1685 1008}%
\special{fp}%
\special{pa 1739 1041}%
\special{pa 1770 1059}%
\special{pa 1775 1061}%
\special{pa 1793 1072}%
\special{fp}%
\special{pa 1848 1105}%
\special{pa 1865 1115}%
\special{pa 1870 1117}%
\special{pa 1902 1136}%
\special{fp}%
\special{pa 1956 1169}%
\special{pa 1960 1170}%
\special{pa 2010 1200}%
\special{fp}%
\special{pn 8}%
\special{pa 2018 1195}%
\special{pa 2040 1182}%
\special{pa 2045 1180}%
\special{pa 2058 1172}%
\special{ip}%
\special{pa 2066 1167}%
\special{pa 2105 1144}%
\special{ip}%
\special{pa 2113 1139}%
\special{pa 2125 1132}%
\special{pa 2130 1130}%
\special{pa 2153 1116}%
\special{ip}%
\special{pa 2161 1111}%
\special{pa 2200 1088}%
\special{ip}%
%
\special{pn 8}%
\special{pa 200 1200}%
\special{pa 2200 1200}%
\special{fp}%
\special{sh 1}%
\special{pa 2200 1200}%
\special{pa 2133 1180}%
\special{pa 2147 1200}%
\special{pa 2133 1220}%
\special{pa 2200 1200}%
\special{fp}%
\put(22.3000,-12.4000){\makebox(0,0)[lb]{{\tiny$x_1$}}}%
\put(2.8000,-13.3000){\makebox(0,0)[lb]{{\tiny$\frac{\pi}{6}$}}}%
\put(20.3000,-13.3000){\makebox(0,0)[lb]{{\tiny$\frac{\pi}{2}$}}}%
%
\special{pn 4}%
\special{pa 400 200}%
\special{pa 400 2200}%
\special{dt 0.045}%
%
\special{pn 4}%
\special{pa 2000 200}%
\special{pa 2000 2200}%
\special{dt 0.045}%
\special{pn 13}%
\special{pa 968 200}%
\special{pa 970 202}%
\special{pa 975 209}%
\special{pa 990 227}%
\special{pa 995 234}%
\special{pa 1010 252}%
\special{pa 1015 259}%
\special{pa 1040 289}%
\special{pa 1045 296}%
\special{pa 1115 380}%
\special{pa 1120 385}%
\special{pa 1145 415}%
\special{pa 1150 420}%
\special{pa 1165 438}%
\special{pa 1170 443}%
\special{pa 1185 461}%
\special{pa 1190 466}%
\special{pa 1200 478}%
\special{pa 1205 483}%
\special{pa 1210 489}%
\special{pa 1215 494}%
\special{pa 1220 500}%
\special{pa 1225 505}%
\special{pa 1230 511}%
\special{pa 1235 516}%
\special{pa 1240 522}%
\special{pa 1245 527}%
\special{pa 1250 533}%
\special{pa 1255 538}%
\special{pa 1260 544}%
\special{pa 1265 549}%
\special{pa 1270 555}%
\special{pa 1275 560}%
\special{pa 1280 566}%
\special{pa 1290 576}%
\special{pa 1295 582}%
\special{pa 1310 597}%
\special{pa 1315 603}%
\special{pa 1330 618}%
\special{pa 1335 624}%
\special{pa 1355 644}%
\special{pa 1360 650}%
\special{pa 1440 730}%
\special{pa 1445 734}%
\special{pa 1465 754}%
\special{pa 1470 758}%
\special{pa 1485 773}%
\special{pa 1490 777}%
\special{pa 1505 792}%
\special{pa 1510 796}%
\special{pa 1515 801}%
\special{pa 1520 805}%
\special{pa 1530 815}%
\special{pa 1535 819}%
\special{pa 1540 824}%
\special{pa 1545 828}%
\special{pa 1550 833}%
\special{pa 1555 837}%
\special{pa 1560 842}%
\special{pa 1565 846}%
\special{pa 1570 851}%
\special{pa 1575 855}%
\special{pa 1580 860}%
\special{pa 1585 864}%
\special{pa 1590 869}%
\special{pa 1595 873}%
\special{pa 1600 878}%
\special{pa 1610 886}%
\special{pa 1615 891}%
\special{pa 1630 903}%
\special{pa 1635 908}%
\special{pa 1650 920}%
\special{pa 1655 925}%
\special{pa 1675 941}%
\special{pa 1680 946}%
\special{pa 1755 1006}%
\special{pa 1760 1009}%
\special{pa 1785 1029}%
\special{pa 1790 1032}%
\special{pa 1805 1044}%
\special{pa 1810 1047}%
\special{pa 1825 1059}%
\special{pa 1830 1062}%
\special{pa 1840 1070}%
\special{pa 1845 1073}%
\special{pa 1850 1077}%
\special{pa 1855 1080}%
\special{pa 1860 1084}%
\special{pa 1865 1087}%
\special{pa 1870 1091}%
\special{pa 1875 1094}%
\special{pa 1880 1098}%
\special{pa 1885 1101}%
\special{pa 1890 1105}%
\special{pa 1895 1108}%
\special{pa 1900 1112}%
\special{pa 1905 1115}%
\special{pa 1910 1119}%
\special{pa 1915 1122}%
\special{pa 1920 1126}%
\special{pa 1930 1132}%
\special{pa 1935 1136}%
\special{pa 1950 1145}%
\special{pa 1955 1149}%
\special{pa 1970 1158}%
\special{pa 1975 1162}%
\special{pa 1995 1174}%
\special{pa 2000 1178}%
\special{fp}%
\special{pn 13}%
\special{pa 2011 1185}%
\special{pa 2039 1202}%
\special{ip}%
\special{pa 2051 1208}%
\special{pa 2079 1225}%
\special{ip}%
\special{pa 2091 1231}%
\special{pa 2105 1240}%
\special{pa 2110 1242}%
\special{pa 2119 1247}%
\special{ip}%
\special{pa 2131 1253}%
\special{pa 2145 1262}%
\special{pa 2150 1264}%
\special{pa 2159 1269}%
\special{ip}%
\special{pa 2171 1275}%
\special{pa 2175 1277}%
\special{pa 2180 1280}%
\special{pa 2185 1282}%
\special{pa 2190 1285}%
\special{pa 2195 1287}%
\special{pa 2200 1290}%
\special{ip}%
\put(12.0000,-4.0000){\makebox(0,0)[lb]{$\widetilde{F}\circ{c_1}^{-1}$}}%
\put(6.3000,-13.3000){\makebox(0,0)[lb]{{\tiny$c_1(t_1)$}}}%
%
\special{pn 8}%
\special{pa 760 550}%
\special{pa 760 1200}%
\special{dt 0.045}%
\end{picture}}%
}
\caption{{\footnotesize The graph of $\widetilde{F}\circ c_1^{-1}$ (Type I)}}
\label{tilde-F-graph1}
\end{figure}
\begin{figure}[H]
\centering
\begin{minipage}{0.45\columnwidth}
\centering
\scalebox{1.0}{\input{tilde-F-graph2.tex}}
\caption{{\footnotesize The graph of $\widetilde{F}\circ c_1^{-1}$ (Type II)}}
\end{minipage}
\begin{minipage}{0.45\columnwidth}
\centering
\scalebox{1.0}{\input{tilde-F-graph3.tex}}
\caption{{\footnotesize The graph of $\widetilde{F}\circ c_1^{-1}$ (Type III)}}
\end{minipage}
\end{figure}
\begin{figure}[H]
\centering
\begin{minipage}{0.45\columnwidth}
\centering
\scalebox{1.0}{\input{tilde-F-graph4.tex}}
\caption{{\footnotesize The graph of $\widetilde{F}\circ c_1^{-1}$ (Type IV)}}
\end{minipage}
\begin{minipage}{0.45\columnwidth}
\centering
\scalebox{1.0}{
{\unitlength 0.1in%
\begin{picture}(20.3000,20.0000)(2.0000,-22.0000)%
\special{pn 8}%
\special{pn 8}%
\special{pa 400 1196}%
\special{pa 405 1190}%
\special{pa 410 1181}%
\special{pa 420 1159}%
\special{pa 425 1146}%
\special{pa 427 1140}%
\special{fp}%
\special{pa 448 1080}%
\special{pa 455 1056}%
\special{pa 460 1039}%
\special{pa 465 1023}%
\special{pa 466 1020}%
\special{fp}%
\special{pa 484 959}%
\special{pa 485 956}%
\special{pa 503 900}%
\special{fp}%
\special{pa 523 839}%
\special{pa 525 832}%
\special{pa 535 804}%
\special{pa 544 781}%
\special{fp}%
\special{pa 569 722}%
\special{pa 570 720}%
\special{pa 575 709}%
\special{pa 585 689}%
\special{pa 595 671}%
\special{pa 597 667}%
\special{fp}%
\special{pa 634 615}%
\special{pa 640 607}%
\special{pa 660 587}%
\special{pa 679 572}%
\special{fp}%
\special{pa 736 546}%
\special{pa 755 542}%
\special{pa 760 542}%
\special{pa 765 541}%
\special{pa 795 541}%
\special{pa 798 542}%
\special{fp}%
\special{pa 861 552}%
\special{pa 865 554}%
\special{pa 870 555}%
\special{pa 875 557}%
\special{pa 880 558}%
\special{pa 890 562}%
\special{pa 895 563}%
\special{pa 919 573}%
\special{fp}%
\special{pa 978 598}%
\special{pa 985 601}%
\special{pa 990 604}%
\special{pa 995 606}%
\special{pa 1000 609}%
\special{pa 1005 611}%
\special{pa 1010 614}%
\special{pa 1015 616}%
\special{pa 1020 619}%
\special{pa 1025 621}%
\special{pa 1033 626}%
\special{fp}%
\special{pa 1089 656}%
\special{pa 1105 665}%
\special{pa 1110 667}%
\special{pa 1135 682}%
\special{pa 1140 684}%
\special{pa 1143 686}%
\special{fp}%
\special{pa 1198 718}%
\special{pa 1252 750}%
\special{fp}%
\special{pa 1307 782}%
\special{pa 1360 814}%
\special{fp}%
\special{pa 1415 847}%
\special{pa 1468 879}%
\special{fp}%
\special{pa 1523 912}%
\special{pa 1576 944}%
\special{fp}%
\special{pa 1631 977}%
\special{pa 1645 985}%
\special{pa 1650 987}%
\special{pa 1685 1008}%
\special{fp}%
\special{pa 1739 1041}%
\special{pa 1770 1059}%
\special{pa 1775 1061}%
\special{pa 1793 1072}%
\special{fp}%
\special{pa 1848 1105}%
\special{pa 1865 1115}%
\special{pa 1870 1117}%
\special{pa 1902 1136}%
\special{fp}%
\special{pa 1956 1169}%
\special{pa 1960 1170}%
\special{pa 2010 1200}%
\special{fp}%
\special{pn 8}%
\special{pa 2018 1195}%
\special{pa 2040 1182}%
\special{pa 2045 1180}%
\special{pa 2058 1172}%
\special{ip}%
\special{pa 2066 1167}%
\special{pa 2105 1144}%
\special{ip}%
\special{pa 2113 1139}%
\special{pa 2125 1132}%
\special{pa 2130 1130}%
\special{pa 2153 1116}%
\special{ip}%
\special{pa 2161 1111}%
\special{pa 2200 1088}%
\special{ip}%
%
\special{pn 8}%
\special{pa 200 1200}%
\special{pa 2200 1200}%
\special{fp}%
\special{sh 1}%
\special{pa 2200 1200}%
\special{pa 2133 1180}%
\special{pa 2147 1200}%
\special{pa 2133 1220}%
\special{pa 2200 1200}%
\special{fp}%
\put(22.3000,-12.4000){\makebox(0,0)[lb]{{\tiny$x_1$}}}%
\put(2.8000,-13.3000){\makebox(0,0)[lb]{{\tiny$\frac{\pi}{6}$}}}%
\put(20.3000,-13.3000){\makebox(0,0)[lb]{{\tiny$\frac{\pi}{2}$}}}%
%
\special{pn 4}%
\special{pa 400 200}%
\special{pa 400 2200}%
\special{dt 0.045}%
%
\special{pn 4}%
\special{pa 2000 200}%
\special{pa 2000 2200}%
\special{dt 0.045}%
\put(12.5000,-21.0000){\makebox(0,0)[lb]{$\widetilde{F}\circ{c_1}^{-1}$}}%
\special{pn 13}%
\special{pa 1106 2200}%
\special{pa 1110 2190}%
\special{pa 1140 2124}%
\special{pa 1145 2114}%
\special{pa 1150 2103}%
\special{pa 1155 2093}%
\special{pa 1160 2082}%
\special{pa 1165 2072}%
\special{pa 1170 2061}%
\special{pa 1175 2051}%
\special{pa 1180 2040}%
\special{pa 1220 1960}%
\special{pa 1225 1951}%
\special{pa 1230 1941}%
\special{pa 1235 1932}%
\special{pa 1240 1922}%
\special{pa 1245 1913}%
\special{pa 1250 1903}%
\special{pa 1255 1894}%
\special{pa 1260 1884}%
\special{pa 1300 1812}%
\special{pa 1305 1804}%
\special{pa 1310 1795}%
\special{pa 1315 1787}%
\special{pa 1320 1778}%
\special{pa 1325 1770}%
\special{pa 1330 1761}%
\special{pa 1335 1753}%
\special{pa 1340 1744}%
\special{pa 1380 1680}%
\special{pa 1385 1673}%
\special{pa 1390 1665}%
\special{pa 1395 1658}%
\special{pa 1400 1650}%
\special{pa 1405 1643}%
\special{pa 1410 1635}%
\special{pa 1415 1628}%
\special{pa 1420 1620}%
\special{pa 1460 1564}%
\special{pa 1465 1558}%
\special{pa 1470 1551}%
\special{pa 1475 1545}%
\special{pa 1480 1538}%
\special{pa 1485 1532}%
\special{pa 1490 1525}%
\special{pa 1495 1519}%
\special{pa 1500 1512}%
\special{pa 1540 1464}%
\special{pa 1545 1459}%
\special{pa 1550 1453}%
\special{pa 1555 1448}%
\special{pa 1560 1442}%
\special{pa 1565 1437}%
\special{pa 1570 1431}%
\special{pa 1575 1426}%
\special{pa 1580 1420}%
\special{pa 1620 1380}%
\special{pa 1625 1376}%
\special{pa 1630 1371}%
\special{pa 1635 1367}%
\special{pa 1640 1362}%
\special{pa 1645 1358}%
\special{pa 1650 1353}%
\special{pa 1655 1349}%
\special{pa 1660 1344}%
\special{pa 1700 1312}%
\special{pa 1705 1309}%
\special{pa 1710 1305}%
\special{pa 1715 1302}%
\special{pa 1720 1298}%
\special{pa 1725 1295}%
\special{pa 1730 1291}%
\special{pa 1735 1288}%
\special{pa 1740 1284}%
\special{pa 1780 1260}%
\special{pa 1785 1258}%
\special{pa 1790 1255}%
\special{pa 1795 1253}%
\special{pa 1800 1250}%
\special{pa 1805 1248}%
\special{pa 1810 1245}%
\special{pa 1815 1243}%
\special{pa 1820 1240}%
\special{pa 1860 1224}%
\special{pa 1865 1223}%
\special{pa 1870 1221}%
\special{pa 1875 1220}%
\special{pa 1880 1218}%
\special{pa 1885 1217}%
\special{pa 1890 1215}%
\special{pa 1895 1214}%
\special{pa 1900 1212}%
\special{pa 1940 1204}%
\special{pa 1945 1204}%
\special{pa 1950 1203}%
\special{pa 1955 1203}%
\special{pa 1960 1202}%
\special{pa 1965 1202}%
\special{pa 1970 1201}%
\special{pa 1975 1201}%
\special{pa 1980 1200}%
\special{pa 2000 1200}%
\special{fp}%
\special{pn 13}%
\special{pa 2015 1200}%
\special{pa 2020 1200}%
\special{pa 2025 1201}%
\special{pa 2030 1201}%
\special{pa 2035 1202}%
\special{pa 2040 1202}%
\special{pa 2045 1203}%
\special{pa 2050 1203}%
\special{pa 2052 1203}%
\special{ip}%
\special{pa 2067 1205}%
\special{pa 2100 1212}%
\special{pa 2103 1213}%
\special{ip}%
\special{pa 2117 1217}%
\special{pa 2120 1218}%
\special{pa 2125 1220}%
\special{pa 2130 1221}%
\special{pa 2135 1223}%
\special{pa 2140 1224}%
\special{pa 2152 1229}%
\special{ip}%
\special{pa 2166 1235}%
\special{pa 2180 1240}%
\special{pa 2185 1243}%
\special{pa 2190 1245}%
\special{pa 2195 1248}%
\special{pa 2200 1250}%
\special{ip}%
%
\special{pn 4}%
\special{pa 760 550}%
\special{pa 760 1200}%
\special{dt 0.045}%
\put(6.3000,-13.3000){\makebox(0,0)[lb]{{\tiny$c_1(t_1)$}}}%
\end{picture}}%
}
\caption{{\footnotesize The graph of $\widetilde{F}\circ c_1^{-1}$ (Type V)}}
\label{tilde-F-graph5}
\end{minipage}
\end{figure}

Hence, from the above classification of the graph of $\widetilde{F}\circ c_1^{-1}$ and $V(x_1,0)=\int_{c_1(t_*)}^{x_1}(\widetilde F\circ c_1^{-1})(x)\,dx+V(c_1(t_*),0)$, we find the shape of the graph of $V(\cdot,0)|_{[\frac{\pi}{6},\frac{\pi}{2}]}$ is as in 
one of Figures \ref{ex-V-graph1}-\ref{ex-V-graph5}.
\vspace{0.5cm}
\begin{figure}[H]
\centering
\scalebox{1.0}{
{\unitlength 0.1in%
\begin{picture}(20.3000,20.0000)(2.0000,-22.0000)%
%
\special{pn 8}%
\special{pa 200 1200}%
\special{pa 2200 1200}%
\special{fp}%
\special{sh 1}%
\special{pa 2200 1200}%
\special{pa 2133 1180}%
\special{pa 2147 1200}%
\special{pa 2133 1220}%
\special{pa 2200 1200}%
\special{fp}%
\put(22.3000,-12.4000){\makebox(0,0)[lb]{{\tiny$x_1$}}}%
\put(2.8000,-13.3000){\makebox(0,0)[lb]{{\tiny$\frac{\pi}{6}$}}}%
\put(20.3000,-13.3000){\makebox(0,0)[lb]{{\tiny$\frac{\pi}{2}$}}}%
%
\special{pn 4}%
\special{pa 400 200}%
\special{pa 400 2200}%
\special{dt 0.045}%
%
\special{pn 4}%
\special{pa 2000 200}%
\special{pa 2000 2200}%
\special{dt 0.045}%
\put(12.0000,-18.0000){\makebox(0,0)[lb]{$V(\cdot,0)|_{[\frac{\pi}{6},\frac{\pi}{2}]}$}}%
\special{pn 13}%
\special{pa 869 2200}%
\special{pa 870 2196}%
\special{pa 900 2112}%
\special{pa 905 2099}%
\special{pa 910 2085}%
\special{pa 915 2072}%
\special{pa 920 2058}%
\special{pa 925 2045}%
\special{pa 930 2031}%
\special{pa 935 2018}%
\special{pa 940 2004}%
\special{pa 980 1900}%
\special{pa 985 1888}%
\special{pa 990 1875}%
\special{pa 995 1863}%
\special{pa 1000 1850}%
\special{pa 1005 1838}%
\special{pa 1010 1825}%
\special{pa 1020 1801}%
\special{pa 1025 1788}%
\special{pa 1055 1716}%
\special{pa 1060 1705}%
\special{pa 1070 1681}%
\special{pa 1075 1670}%
\special{pa 1080 1658}%
\special{pa 1085 1647}%
\special{pa 1090 1635}%
\special{pa 1095 1624}%
\special{pa 1100 1612}%
\special{pa 1140 1524}%
\special{pa 1145 1514}%
\special{pa 1150 1503}%
\special{pa 1155 1493}%
\special{pa 1160 1482}%
\special{pa 1165 1472}%
\special{pa 1170 1461}%
\special{pa 1175 1451}%
\special{pa 1180 1440}%
\special{pa 1220 1360}%
\special{pa 1225 1351}%
\special{pa 1230 1341}%
\special{pa 1235 1332}%
\special{pa 1240 1322}%
\special{pa 1245 1313}%
\special{pa 1250 1303}%
\special{pa 1255 1294}%
\special{pa 1260 1284}%
\special{pa 1300 1212}%
\special{pa 1305 1204}%
\special{pa 1310 1195}%
\special{pa 1315 1187}%
\special{pa 1320 1178}%
\special{pa 1325 1170}%
\special{pa 1330 1161}%
\special{pa 1335 1153}%
\special{pa 1340 1144}%
\special{pa 1380 1080}%
\special{pa 1385 1073}%
\special{pa 1390 1065}%
\special{pa 1395 1058}%
\special{pa 1400 1050}%
\special{pa 1405 1043}%
\special{pa 1410 1035}%
\special{pa 1415 1028}%
\special{pa 1420 1020}%
\special{pa 1460 964}%
\special{pa 1465 958}%
\special{pa 1470 951}%
\special{pa 1475 945}%
\special{pa 1480 938}%
\special{pa 1485 932}%
\special{pa 1490 925}%
\special{pa 1495 919}%
\special{pa 1500 912}%
\special{pa 1540 864}%
\special{pa 1545 859}%
\special{pa 1550 853}%
\special{pa 1555 848}%
\special{pa 1560 842}%
\special{pa 1565 837}%
\special{pa 1570 831}%
\special{pa 1575 826}%
\special{pa 1580 820}%
\special{pa 1620 780}%
\special{pa 1625 776}%
\special{pa 1630 771}%
\special{pa 1635 767}%
\special{pa 1640 762}%
\special{pa 1645 758}%
\special{pa 1650 753}%
\special{pa 1655 749}%
\special{pa 1660 744}%
\special{pa 1700 712}%
\special{pa 1705 709}%
\special{pa 1710 705}%
\special{pa 1715 702}%
\special{pa 1720 698}%
\special{pa 1725 695}%
\special{pa 1730 691}%
\special{pa 1735 688}%
\special{pa 1740 684}%
\special{pa 1780 660}%
\special{pa 1785 658}%
\special{pa 1790 655}%
\special{pa 1795 653}%
\special{pa 1800 650}%
\special{pa 1805 648}%
\special{pa 1810 645}%
\special{pa 1815 643}%
\special{pa 1820 640}%
\special{pa 1860 624}%
\special{pa 1865 623}%
\special{pa 1870 621}%
\special{pa 1875 620}%
\special{pa 1880 618}%
\special{pa 1885 617}%
\special{pa 1890 615}%
\special{pa 1895 614}%
\special{pa 1900 612}%
\special{pa 1940 604}%
\special{pa 1945 604}%
\special{pa 1950 603}%
\special{pa 1955 603}%
\special{pa 1960 602}%
\special{pa 1965 602}%
\special{pa 1970 601}%
\special{pa 1975 601}%
\special{pa 1980 600}%
\special{pa 2000 600}%
\special{fp}%
\special{pn 13}%
\special{pa 2015 600}%
\special{pa 2020 600}%
\special{pa 2025 601}%
\special{pa 2030 601}%
\special{pa 2035 602}%
\special{pa 2040 602}%
\special{pa 2045 603}%
\special{pa 2050 603}%
\special{pa 2052 603}%
\special{ip}%
\special{pa 2067 605}%
\special{pa 2100 612}%
\special{pa 2103 613}%
\special{ip}%
\special{pa 2117 617}%
\special{pa 2120 618}%
\special{pa 2125 620}%
\special{pa 2130 621}%
\special{pa 2135 623}%
\special{pa 2140 624}%
\special{pa 2152 629}%
\special{ip}%
\special{pa 2166 635}%
\special{pa 2180 640}%
\special{pa 2185 643}%
\special{pa 2190 645}%
\special{pa 2195 648}%
\special{pa 2200 650}%
\special{ip}%
%
\special{pn 4}%
\special{pa 800 1200}%
\special{pa 800 2200}%
\special{dt 0.045}%
\end{picture}}%
}
\caption{{\footnotesize The graph of $V(\cdot,0)|_{[\frac{\pi}{6},\frac{\pi}{2}]}$(Type I)}}
\label{ex-V-graph1}
\end{figure}
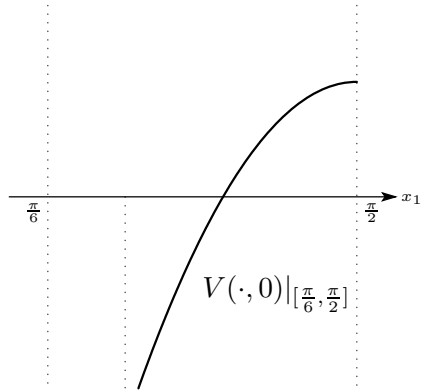
\begin{figure}[H]
\centering
\begin{minipage}{0.45\columnwidth}
\centering
\scalebox{1.0}{\input{ex-V-graph2.tex}}
\caption{{\footnotesize The graph of $V(\cdot,0)|_{[\frac{\pi}{6},\frac{\pi}{2}]}$(Type II)}}
\end{minipage}
\hspace{0.3cm}
\begin{minipage}{0.45\columnwidth}
\centering
\scalebox{1.0}{
{\unitlength 0.1in%
\begin{picture}(20.3000,20.0000)(2.0000,-22.0000)%
%
\special{pn 8}%
\special{pa 200 1200}%
\special{pa 2200 1200}%
\special{fp}%
\special{sh 1}%
\special{pa 2200 1200}%
\special{pa 2133 1180}%
\special{pa 2147 1200}%
\special{pa 2133 1220}%
\special{pa 2200 1200}%
\special{fp}%
\put(22.3000,-12.4000){\makebox(0,0)[lb]{{\tiny$x_1$}}}%
\put(2.8000,-13.3000){\makebox(0,0)[lb]{{\tiny$\frac{\pi}{6}$}}}%
\put(20.4000,-13.3000){\makebox(0,0)[lb]{{\tiny$\frac{\pi}{2}$}}}%
%
\special{pn 4}%
\special{pa 400 200}%
\special{pa 400 2200}%
\special{dt 0.045}%
%
\special{pn 4}%
\special{pa 2000 200}%
\special{pa 2000 2200}%
\special{dt 0.045}%
\put(5.7000,-8.7000){\makebox(0,0)[lb]{$V(\cdot,0)|_{[\frac{\pi}{6},\frac{\pi}{2}]}$}}%
\special{pn 13}%
\special{pn 13}%
\special{pa 200 1754}%
\special{pa 205 1757}%
\special{pa 214 1760}%
\special{ip}%
\special{pa 248 1773}%
\special{pa 255 1776}%
\special{pa 260 1777}%
\special{pa 262 1778}%
\special{ip}%
\special{pa 298 1788}%
\special{pa 300 1788}%
\special{pa 305 1790}%
\special{pa 312 1791}%
\special{ip}%
\special{pa 348 1797}%
\special{pa 350 1797}%
\special{pa 355 1798}%
\special{pa 360 1798}%
\special{pa 363 1799}%
\special{ip}%
\special{ip}%
\special{pa 400 1800}%
\special{pa 420 1800}%
\special{pa 425 1799}%
\special{pa 435 1799}%
\special{pa 440 1798}%
\special{pa 445 1798}%
\special{pa 450 1797}%
\special{pa 455 1797}%
\special{pa 470 1794}%
\special{pa 475 1794}%
\special{pa 495 1790}%
\special{pa 500 1788}%
\special{pa 515 1785}%
\special{pa 520 1783}%
\special{pa 530 1781}%
\special{pa 540 1777}%
\special{pa 545 1776}%
\special{pa 555 1772}%
\special{pa 560 1771}%
\special{pa 595 1757}%
\special{pa 600 1754}%
\special{pa 610 1750}%
\special{pa 615 1747}%
\special{pa 620 1745}%
\special{pa 625 1742}%
\special{pa 630 1740}%
\special{pa 635 1737}%
\special{pa 640 1735}%
\special{pa 660 1723}%
\special{pa 665 1721}%
\special{pa 680 1712}%
\special{pa 685 1708}%
\special{pa 705 1696}%
\special{pa 710 1692}%
\special{pa 715 1689}%
\special{pa 720 1685}%
\special{pa 725 1682}%
\special{pa 730 1678}%
\special{pa 735 1675}%
\special{pa 740 1671}%
\special{pa 745 1668}%
\special{pa 765 1652}%
\special{pa 770 1649}%
\special{pa 790 1633}%
\special{pa 795 1628}%
\special{pa 815 1612}%
\special{pa 820 1607}%
\special{pa 830 1599}%
\special{pa 835 1594}%
\special{pa 840 1590}%
\special{pa 845 1585}%
\special{pa 850 1581}%
\special{pa 860 1571}%
\special{pa 865 1567}%
\special{pa 875 1557}%
\special{pa 880 1553}%
\special{pa 910 1523}%
\special{pa 915 1519}%
\special{pa 920 1513}%
\special{pa 955 1478}%
\special{pa 960 1472}%
\special{pa 975 1457}%
\special{pa 980 1451}%
\special{pa 985 1446}%
\special{pa 990 1440}%
\special{pa 1000 1430}%
\special{pa 1005 1424}%
\special{pa 1010 1419}%
\special{pa 1015 1413}%
\special{pa 1020 1408}%
\special{pa 1025 1402}%
\special{pa 1030 1397}%
\special{pa 1040 1385}%
\special{pa 1045 1380}%
\special{pa 1050 1374}%
\special{pa 1055 1369}%
\special{pa 1070 1351}%
\special{pa 1075 1346}%
\special{pa 1085 1334}%
\special{pa 1090 1329}%
\special{pa 1110 1305}%
\special{pa 1115 1300}%
\special{pa 1135 1276}%
\special{pa 1140 1271}%
\special{pa 1175 1229}%
\special{pa 1180 1224}%
\special{pa 1220 1176}%
\special{pa 1225 1171}%
\special{pa 1260 1129}%
\special{pa 1265 1124}%
\special{pa 1285 1100}%
\special{pa 1290 1095}%
\special{pa 1310 1071}%
\special{pa 1315 1066}%
\special{pa 1325 1054}%
\special{pa 1330 1049}%
\special{pa 1345 1031}%
\special{pa 1350 1026}%
\special{pa 1355 1020}%
\special{pa 1360 1015}%
\special{pa 1370 1003}%
\special{pa 1375 998}%
\special{pa 1380 992}%
\special{pa 1385 987}%
\special{pa 1390 981}%
\special{pa 1395 976}%
\special{pa 1400 970}%
\special{pa 1410 960}%
\special{pa 1415 954}%
\special{pa 1420 949}%
\special{pa 1425 943}%
\special{pa 1440 928}%
\special{pa 1445 922}%
\special{pa 1480 887}%
\special{pa 1485 881}%
\special{pa 1490 877}%
\special{pa 1520 847}%
\special{pa 1525 843}%
\special{pa 1535 833}%
\special{pa 1540 829}%
\special{pa 1550 819}%
\special{pa 1555 815}%
\special{pa 1560 810}%
\special{pa 1565 806}%
\special{pa 1570 801}%
\special{pa 1580 793}%
\special{pa 1585 788}%
\special{pa 1605 772}%
\special{pa 1610 767}%
\special{pa 1630 751}%
\special{pa 1635 748}%
\special{pa 1655 732}%
\special{pa 1660 729}%
\special{pa 1665 725}%
\special{pa 1670 722}%
\special{pa 1675 718}%
\special{pa 1680 715}%
\special{pa 1685 711}%
\special{pa 1690 708}%
\special{pa 1695 704}%
\special{pa 1715 692}%
\special{pa 1720 688}%
\special{pa 1735 679}%
\special{pa 1740 677}%
\special{pa 1760 665}%
\special{pa 1765 663}%
\special{pa 1770 660}%
\special{pa 1775 658}%
\special{pa 1780 655}%
\special{pa 1785 653}%
\special{pa 1790 650}%
\special{pa 1800 646}%
\special{pa 1805 643}%
\special{pa 1840 629}%
\special{pa 1845 628}%
\special{pa 1855 624}%
\special{pa 1860 623}%
\special{pa 1870 619}%
\special{pa 1880 617}%
\special{pa 1885 615}%
\special{pa 1900 612}%
\special{pa 1905 610}%
\special{pa 1925 606}%
\special{pa 1930 606}%
\special{pa 1945 603}%
\special{pa 1950 603}%
\special{pa 1955 602}%
\special{pa 1960 602}%
\special{pa 1965 601}%
\special{pa 1975 601}%
\special{pa 1980 600}%
\special{pa 2000 600}%
\special{fp}%
\special{pn 13}%
\special{pa 2015 600}%
\special{pa 2020 600}%
\special{pa 2025 601}%
\special{pa 2035 601}%
\special{pa 2040 602}%
\special{pa 2045 602}%
\special{pa 2050 603}%
\special{pa 2052 603}%
\special{ip}%
\special{pa 2066 605}%
\special{pa 2070 606}%
\special{pa 2075 606}%
\special{pa 2095 610}%
\special{pa 2100 612}%
\special{pa 2102 612}%
\special{ip}%
\special{pa 2117 616}%
\special{pa 2120 617}%
\special{pa 2130 619}%
\special{pa 2140 623}%
\special{pa 2145 624}%
\special{pa 2152 627}%
\special{ip}%
\special{pa 2166 631}%
\special{pa 2195 643}%
\special{pa 2200 646}%
\special{ip}%
\end{picture}}%
}
\caption{{\footnotesize The graph of $V(\cdot,0)|_{[\frac{\pi}{6},\frac{\pi}{2}]}$(Type III)}}
\end{minipage}
\end{figure}
\begin{figure}[H]
\centering
\begin{minipage}{0.45\columnwidth}
\centering
\scalebox{1.0}{\input{ex-V-graph4.tex}}
\caption{{\footnotesize The graph of $V(\cdot,0)|_{[\frac{\pi}{6},\frac{\pi}{2}]}$(Type IV)}}
\end{minipage}
\hspace{0.3cm}
\begin{minipage}{0.45\columnwidth}
\centering
\scalebox{1.0}{
{\unitlength 0.1in%
\begin{picture}(20.3000,20.0000)(2.0000,-22.0000)%
%
\special{pn 8}%
\special{pa 200 1200}%
\special{pa 2200 1200}%
\special{fp}%
\special{sh 1}%
\special{pa 2200 1200}%
\special{pa 2133 1180}%
\special{pa 2147 1200}%
\special{pa 2133 1220}%
\special{pa 2200 1200}%
\special{fp}%
\put(22.3000,-12.4000){\makebox(0,0)[lb]{{\tiny$x_1$}}}%
\put(2.8000,-13.3000){\makebox(0,0)[lb]{{\tiny$\frac{\pi}{6}$}}}%
\put(20.3000,-13.3000){\makebox(0,0)[lb]{{\tiny$\frac{\pi}{2}$}}}%
%
\special{pn 4}%
\special{pa 400 200}%
\special{pa 400 2200}%
\special{dt 0.045}%
%
\special{pn 4}%
\special{pa 2000 200}%
\special{pa 2000 2200}%
\special{dt 0.045}%
\put(10.1000,-4.5000){\makebox(0,0)[lb]{$V(\cdot,0)|_{[\frac{\pi}{6},\frac{\pi}{2}]}$}}%
\special{pn 13}%
\special{pa 869 200}%
\special{pa 870 204}%
\special{pa 900 288}%
\special{pa 905 301}%
\special{pa 910 315}%
\special{pa 915 328}%
\special{pa 920 342}%
\special{pa 925 355}%
\special{pa 930 369}%
\special{pa 935 382}%
\special{pa 940 396}%
\special{pa 980 500}%
\special{pa 985 512}%
\special{pa 990 525}%
\special{pa 995 537}%
\special{pa 1000 550}%
\special{pa 1005 562}%
\special{pa 1010 575}%
\special{pa 1020 599}%
\special{pa 1025 612}%
\special{pa 1055 684}%
\special{pa 1060 695}%
\special{pa 1070 719}%
\special{pa 1075 730}%
\special{pa 1080 742}%
\special{pa 1085 753}%
\special{pa 1090 765}%
\special{pa 1095 776}%
\special{pa 1100 788}%
\special{pa 1140 876}%
\special{pa 1145 886}%
\special{pa 1150 897}%
\special{pa 1155 907}%
\special{pa 1160 918}%
\special{pa 1165 928}%
\special{pa 1170 939}%
\special{pa 1175 949}%
\special{pa 1180 960}%
\special{pa 1220 1040}%
\special{pa 1225 1049}%
\special{pa 1230 1059}%
\special{pa 1235 1068}%
\special{pa 1240 1078}%
\special{pa 1245 1087}%
\special{pa 1250 1097}%
\special{pa 1255 1106}%
\special{pa 1260 1116}%
\special{pa 1300 1188}%
\special{pa 1305 1196}%
\special{pa 1310 1205}%
\special{pa 1315 1213}%
\special{pa 1320 1222}%
\special{pa 1325 1230}%
\special{pa 1330 1239}%
\special{pa 1335 1247}%
\special{pa 1340 1256}%
\special{pa 1380 1320}%
\special{pa 1385 1327}%
\special{pa 1390 1335}%
\special{pa 1395 1342}%
\special{pa 1400 1350}%
\special{pa 1405 1357}%
\special{pa 1410 1365}%
\special{pa 1415 1372}%
\special{pa 1420 1380}%
\special{pa 1460 1436}%
\special{pa 1465 1442}%
\special{pa 1470 1449}%
\special{pa 1475 1455}%
\special{pa 1480 1462}%
\special{pa 1485 1468}%
\special{pa 1490 1475}%
\special{pa 1495 1481}%
\special{pa 1500 1488}%
\special{pa 1540 1536}%
\special{pa 1545 1541}%
\special{pa 1550 1547}%
\special{pa 1555 1552}%
\special{pa 1560 1558}%
\special{pa 1565 1563}%
\special{pa 1570 1569}%
\special{pa 1575 1574}%
\special{pa 1580 1580}%
\special{pa 1620 1620}%
\special{pa 1625 1624}%
\special{pa 1630 1629}%
\special{pa 1635 1633}%
\special{pa 1640 1638}%
\special{pa 1645 1642}%
\special{pa 1650 1647}%
\special{pa 1655 1651}%
\special{pa 1660 1656}%
\special{pa 1700 1688}%
\special{pa 1705 1691}%
\special{pa 1710 1695}%
\special{pa 1715 1698}%
\special{pa 1720 1702}%
\special{pa 1725 1705}%
\special{pa 1730 1709}%
\special{pa 1735 1712}%
\special{pa 1740 1716}%
\special{pa 1780 1740}%
\special{pa 1785 1742}%
\special{pa 1790 1745}%
\special{pa 1795 1747}%
\special{pa 1800 1750}%
\special{pa 1805 1752}%
\special{pa 1810 1755}%
\special{pa 1815 1757}%
\special{pa 1820 1760}%
\special{pa 1860 1776}%
\special{pa 1865 1777}%
\special{pa 1870 1779}%
\special{pa 1875 1780}%
\special{pa 1880 1782}%
\special{pa 1885 1783}%
\special{pa 1890 1785}%
\special{pa 1895 1786}%
\special{pa 1900 1788}%
\special{pa 1940 1796}%
\special{pa 1945 1796}%
\special{pa 1950 1797}%
\special{pa 1955 1797}%
\special{pa 1960 1798}%
\special{pa 1965 1798}%
\special{pa 1970 1799}%
\special{pa 1975 1799}%
\special{pa 1980 1800}%
\special{pa 2000 1800}%
\special{fp}%
\special{pn 13}%
\special{pa 2015 1800}%
\special{pa 2020 1800}%
\special{pa 2025 1799}%
\special{pa 2030 1799}%
\special{pa 2035 1798}%
\special{pa 2040 1798}%
\special{pa 2045 1797}%
\special{pa 2050 1797}%
\special{pa 2052 1797}%
\special{ip}%
\special{pa 2067 1795}%
\special{pa 2100 1788}%
\special{pa 2103 1787}%
\special{ip}%
\special{pa 2117 1783}%
\special{pa 2120 1782}%
\special{pa 2125 1780}%
\special{pa 2130 1779}%
\special{pa 2135 1777}%
\special{pa 2140 1776}%
\special{pa 2152 1771}%
\special{ip}%
\special{pa 2166 1765}%
\special{pa 2180 1760}%
\special{pa 2185 1757}%
\special{pa 2190 1755}%
\special{pa 2195 1752}%
\special{pa 2200 1750}%
\special{ip}%
%
\special{pn 4}%
\special{pa 790 200}%
\special{pa 790 1200}%
\special{dt 0.045}%
\end{picture}}%
}
\caption{{\footnotesize The graph of $V(\cdot,0)|_{[\frac{\pi}{6},\frac{\pi}{2}]}$(Type V)}}
\label{ex-V-graph5}
\end{minipage}
\end{figure}

\vspace{1.0cm}
\begin{flushleft}
\hspace{7.53307cm}Tomoki Fujii\\
\hspace{7.53307cm}Department of Mathematics\\
\hspace{7.53307cm}Graduate School of Science\\
\hspace{7.53307cm}Tokyo University of Sience\\
\hspace{7.53307cm}1-3 Kagurazaka, Shinjuku-ku, Tokyo, 162-8601\\
\hspace{7.53307cm}Japan\\
\hspace{7.53307cm}E-mail:1121703@ed.tus.ac.jp
\end{flushleft}
\begin{flushleft}
\hspace{7.53307cm}Naoyuki Koike\\
\hspace{7.53307cm}Department of Mathematics\\
\hspace{7.53307cm}Faculty of Science\\
\hspace{7.53307cm}Tokyo University of Sience\\
\hspace{7.53307cm}1-3 Kagurazaka, Shinjuku-ku, Tokyo, 162-8601\\
\hspace{7.53307cm}Japan\\
\hspace{7.53307cm}E-mail:koike@rs.tus.ac.jp
\end{flushleft}
\end{document}